\documentclass[10pt]{article}
  \usepackage{amsfonts} \usepackage{latexsym}
   \RequirePackage{amsbsy} 
   \RequirePackage{amsopn} 
   \RequirePackage{amsmath}
    \RequirePackage{amssymb}

    \DeclareMathOperator{\Spec}   {Spec}
   \DeclareMathOperator{\Max}    {Max}

   \DeclareMathOperator{\calM}   {\mathcal M}
   \DeclareMathOperator{\calW}   {\mathcal W}

    \DeclareMathOperator{\calR} {\mathcal R}
  \usepackage[T1]{fontenc}

  \usepackage{color}
   \newtheorem{thee}{Theorem}[section]
   \newtheorem{coor}[thee]{Corollary}
   \newtheorem{leem}[thee]{Lemma}
   
   \newtheorem{prro}[thee]{Proposition}
   
   \newtheorem{exxe}[thee]{Example}
   \newtheorem{reem}[thee]{Remark}
  \DeclareMathOperator{\oF}       {\boldsymbol{\overline{F}}}%
  \DeclareMathOperator{\F}        {\boldsymbol{F}}%

   \newcommand{\balf}
   {\renewcommand{\theenumi}{(\alph{enumi})}
   \renewcommand{\labelenumi}{\theenumi}
                        \begin{enumerate}}
  \newcommand{\ealf}   {\end{enumerate}
                        \renewcommand{\theenumi}{\arabic{enumi}}
                        \renewcommand{\labelenumi}{\theenumi.}}
  \newcommand{\bara}   {\renewcommand{\theenumi}{(\arabic{enumi})}
                        \renewcommand{\labelenumi}{\theenumi}
                        \begin{enumerate} }
  \newcommand{\eara}   {\end{enumerate}
                        \renewcommand{\theenumi}{\arabic{enumi}}
                        \renewcommand{\labelenumi}{\theenumi.}}

   \newcommand{\brom}   {\renewcommand{\theenumi}{(\roman{enumi})}
                        \renewcommand{\labelenumi}{\theenumi}
                        \begin{enumerate} }
  \newcommand{\erom}   {\end{enumerate}
                        \renewcommand{\theenumi}{\arabic{enumi}}
                        \renewcommand{\labelenumi}{\theenumi.}}

\begin{document}

   \title {STAR OPERATIONS AND PULLBACKS}

   \author{ {\Large Marco Fontana\thanks{Supported in part by research
   grants MIUR 2001/2002 and 2003/2004.}\ , \;\; Mi Hee Park
   } \\
   \vspace{3pt}\\
  Dipartimento di Matematica \\
  Universit\`a degli Studi Roma Tre \\
  {\small Largo San Leonardo Murialdo, 1 }\\
  {\small 00146 Roma, Italy} \\
   {\small fontana@mat.uniroma3.it} \\
   \vspace{3pt}\\
  School of Mathematics \\
  Korea Institute for Advanced Study \\
  {\small Seoul 130-722, Korea }\\
  {\small mhpark@kias.re.kr}\\
  }

 \date{ } \vskip- 10pt

  \maketitle

   \begin{abstract} In this paper we study the star operations on
   a pullback of integral domains.  In particular, we characterize the
star
   operations of a domain arising from a pullback of ``a general type''
by
   introducing new techniques for ``projecting'' and ``lifting'' star
   operations under surjective homomorphisms of integral domains.  We
study
   the transfer in a pullback (or with respect to a surjective
homomorphism) of
   some relevant classes or distinguished properties of star operations
   such as\ $v-,\ t-,\ w-,\ b-,\ d-,$\ finite type, e.a.b., stable, and
spectral
   operations.\ We apply part of the theory developed here to give a
   complete positive answer to a problem posed by D. F. Anderson in 1992
   concerning the star operations on the ``$D+M$'' constructions.
    \end{abstract}

   \bigskip

  \section{Introduction and preliminary results}

    The theory of ideal systems and star operations was developed by W.
    Krull, H. Pr\"{u}fer, and E. Noether around 1930, and is a
powerful tool
    for characterizing se\-ve\-ral re\-le\-vant classes of integral
domains, for
    studying their mutual relations and for introducing the Kronecker
    function rings in a very general ring-theoretical setting.  A modern
    treatment of various aspects of this theory can be found in the
volumes by
    P. Jaffard \cite{J}, O. Zariski and P. Samuel \cite[Appendix 4]{ZS},
R.
    Gilmer \cite{G}, M.D. Larsen and P.J. McCarthy \cite{LMc}, and F.
    Halter-Koch \cite {HK2}.

    Pullbacks were considered in \cite{F} for providing an appropriate
unified
    setting for several important ``composite--type'' constructions
    introduced in va\-ri\-ous contexts of commutative ring theory in
order
    to construct examples and counter-examples with different
pathologies:\
    for instance, Seidenberg's constructions for (polynomial)
dimensional
    sequences \cite{S:1953}, Nagata's composition of va\-lua\-tion
domains
    and ``$K+J(R)$'' constructions \cite[page 35 and Appendix A1,
Example
    2]{N:1964}, Akiba's AV-domains or Dobbs' divided domains
    \cite{Akiba:1967, Dobbs:1976}, Gilmer's ``$D+M$'' constructions
    \cite{G}, Traverso's glueings for a constructive approach to the
    seminormalization \cite{T:1970}, Vasconcelos' umbrella rings and
    Greenberg's F-domains \cite{Vasconcelos:1972, Greenberg:1978},
    Boisen-Sheldon's CPI-extensions \cite{BS:1977}, Hedstrom-Houston's
    pseudo-valuation domains \cite{HH}, ``$D+XD_{S}[X]$'' rings and
more generally, the ``$A+XB[X]$'' rings considered by many authors (see
the
    recent excellent survey papers by T. Lucas \cite{L:2000} and M.
    Zafrullah \cite{Z:2003}, which contain  ample and updated
bibliographies
    on this subject).

    It was natural at this stage of knowledge to investigate the
behaviour of
    the star operations in a general pullback setting and with respect
to
    surjective homomorphisms of integral domains, after various
different
    results concerning distinguished star operations (like the $v$--,
the
    $t$-- or the $w$-- operation) and particular ``composite--type''
    constructions were obtained by different authors (cf.  for instance
    \cite{BG:1973}, \cite{BR:1976}, \cite{AA:1988}, \cite{DHLZ:1989},
    \cite{Kang:1989},  \cite{AA:1990}, \cite{AAZ:1991}, \cite{AHZ:1993},
\cite{FG:1996},
    \cite{Gabelli/Houston},
    \cite{AElK:1999}, \cite{P},  \cite{Mimouni:2003}, and the survey
papers
    \cite{A:2000} and \cite{GH:2000}).

    The present work was stimulated by the papers by D.D.
    Anderson and D.F. Anderson on  star operations, and more precisely,
    by the study initiated by D.F. Anderson concerning the star
operations
    on the ``$D+M$'' constructions \cite{A:1992}.

    In Section 2, after introducing an operation of ``glueing'' of star
operations
    in a pullback of integral domains, we will characterize the star
    operations of a domain arising from a pullback of ``a general
type''.  For
    this purpose we will introduce new techniques for ``projecting'' and
    ``lifting'' star operations under surjective homomorphisms of
integral
    domains.  Section 3 is devoted to the study of the transfer in a
pullback
    (or with respect to a surjective homomorphism) of some relevant
    properties or classes of star operations such as\ $v$--,\ $t$--,\
$w$--,\
    $b$--,\ $d$--,\ finite type, e.a.b., stable, and spectral operations.

    We will apply part of the theory developed here to give a complete
    positive answer to a problem posed by D. F. Anderson in 1992
    \cite{A:1992} concerning the star operations on the ``$D+M$''
    constructions.

    \bigskip

   Let $D$ be an integral domain with quotient field $L$.  Let
   $\boldsymbol{\overline{F}}(D)$ denote the set of all nonzero
   $D$-submodules of $L$ and let $\boldsymbol{F}(D)$ be the set of all
   nonzero fractional ideals of $D$, i.e., all $E \in
   \boldsymbol{\overline{F}}(D) $ such that there exists a nonzero $d
\in
   D$ with $dE \subseteq D$.  Let $\boldsymbol{f}(D)$ be the set of all
   nonzero finitely generated $D$-submodules of $L$.  Then obviously,
$\,
   \boldsymbol{f}(D) \subseteq \boldsymbol{F}(D) \subseteq
   \boldsymbol{\overline{F}}(D) \, .$

   For each pair of nonzero fractional ideals $E, F$ of $D$, we denote as usual
by
  $(E:_{L}F)$ the fractional ideal of $D$ given by $\{ y\in L\mid \, yF
  \subseteq E\}$; in particular, for each nonzero  fractional ideal $I$
  of $D$, we set $I^{-1} := (D:_{L} I)$.
  \medskip

  We recall that a mapping
  $ \star : \boldsymbol{\overline{F}}(D) \rightarrow
  \boldsymbol{\overline{F}}(D) \,, \; E \mapsto E^{\star} $,
   is called a \it semistar operation on $D$ \rm if the following
   properties hold
   for all $0\neq x \in L$ and $E,F \in \boldsymbol{\overline{F}}(D)$:

  \hspace*{10pt}$(\star_1)$ \;\;  $(xE)^{\star} = xE^{\star}\,; $

  \hspace*{10pt}$(\star_2)$ \;\;  $E \subseteq F
  \Rightarrow  E^{\star} \subseteq F^{\star}\,; $

  \hspace*{10pt}$(\star_3)$\;\;  $E \subseteq E^{\star}$ and
  $E^{\star} =
  (E^{\star})^{\star} =: E^{\star \star} $

  \noindent (cf.  for instance \cite{OM1}, \cite{OM2},
  \cite{MSu}, \cite{MS}, \cite{FH} and \cite{FL1}).

  \begin{exxe}\label{ex:1.1}
  \bf (a) \rm If $\star$ is a semistar operation on $D$ such that
$D^\star =
  D\,$, then the map (still denoted by) $\star : \boldsymbol{F}(D)
  \rightarrow \boldsymbol{F}(D)\,$, \ $E \mapsto E^\star\,,$ is called a
  \it star operation on \rm $D\,.$ Recall \cite[(32.1)]{G} that a star
  operation $\star$ satisfies the properties $(\star_2)\,, (\star_3)\,$
  for all $E, F \in \boldsymbol{F}(D)\,;$ moreover, for each $0\neq x
\in
  L\,$  and $E \in \boldsymbol{F}(D)\,,$ a star
  operation $\star$ satisfies the following:

  \hspace*{10pt} $(\star\star_1)$ $(xD)^\star = xD\,; \; \; (xE)^{\star}
= xE^{\star}\,.$

  \noindent A semistar operation on $D$ such that $D \subsetneq
  D^\star$ is called a \it proper semistar o\-pe\-ration on $D$.

  \bf (b) \rm The \it trivial semistar operation $e_{D}$ on $\,D\,$ \rm
  (simply denoted by $e$) is the semistar operation constant onto
$\,L\,$,
  i.e., the semistar operation defined by $\, E^{e_{D}} := L\,$ for each
  $\, E \in \boldsymbol{\overline{F}}(D)\,$.  \ Note that $\,\star\,$ is
  the trivial semistar operation on $\,D\,$ if and only if $\, D^{\star}
=
  L\,$.

  \bf (c) \rm Another trivial semistar (in fact, star) operation is the
\it
  identity star operation $d_{D}$ on $\,D\,$ \rm (simply denoted by $d$)
  defined by $\, E^{d_{D}}: = E\,$ for each $\, E \in
  \boldsymbol{\overline{F}}(D)\,$.

   \bf (d) \rm For each $E \in \boldsymbol{\overline{F}}(D)$, set
  $E^{\star_f} := \cup \{F^{\star} \;|\;  \, F\subseteq E, \; F \in
  \boldsymbol{f}(D) \}\,.  $ \ Then $\star_f$ is also a semistar
operation
  on $D$, which is called \it the semistar operation of finite type
  associated to $\star\,$.  \rm Obviously, $F^{\star} = F^{\star_{f}}$
for
  each $F \in \boldsymbol{f}(D)\, ;$ \ moreover, if $\star$ is a star
  operation, then ${\star_{f}}$ is also a star operation.  If $\star =
  \star_f$, then the semistar [respectively, the star] operation $\star$
  is called a \it semistar \rm [respectively, \it star\rm ] \it
operation
  of finite type \rm \cite[Example 2.5 (4)]{FL1}.

  Note that, in general, $\star_f \leq \star\,,$ i.e., $\,
  E^{\star_f} \subseteq E^{\star}\, $ for each $E \in
  \boldsymbol{\overline{F}}(D)$. Thus, in particular, if $E =
  E^{\star}$, then $E = E^{\star_f}$.  Note also that $\star_f =
  (\star_f)_{f}$.

  There are several examples of nontrivial semistar operations of finite
type;
  the best known is probably the $t$--operation.  Indeed, we start from
  the \it $v_{D}$ star operation \rm on an integral domain $D$ (simply
  denoted by $v$), which is defined by

  \centerline{ $E^{v_{D}} := (E^{-1})^{-1}=(D:_{L}(D:_{L}E)) $}

  \noindent for any $E\in\F(D)$, and we set  $t_{D}:=(v_{D})_f$ (or
simply, $t
  =v_{f}$).

  Other relevant examples of semistar operations of finite type will be
  constructed later.
  \end{exxe}

  A semistar operation $\star$ on $D$ is called an
  \it e.a.b.\rm \ (\it endlich
  arithmetisch brauchbar\rm ) [respectively,
  \it a.b. \rm (\it arithmetisch brauchbar\rm )] \it semistar operation
\rm  if

  \centerline{$(EF)^\star \subseteq (EG)^\star \Rightarrow F^\star
\subseteq
  G^\star$}

  \noindent for each $E \in {\boldsymbol f}(D)$ and all $F,G \in
{\boldsymbol
  f}(D)$ [respectively, $F,G \in {\overline{\boldsymbol F}}(D)]$
\cite[Definition 2.3 and Lemma 2.7]{FL1}.

  If $\star$ is a star operation on $D$, then the definition of
  e.a.b. [respectively,
   a.b.] operation is analogous (for an a.b. star operation,\ $F, G$\ are
taken in
   ${\boldsymbol F}(D)$).

  \medskip
  \begin{exxe} \label{ex:1.2} \rm Let $\iota:R \hookrightarrow T$ be an
embedding of
  integral domains with the same field of quotients $K$ and let $\ast$
be
  a semistar operation on $R$.  Define $\ast_{\iota} :
  \overline{\boldsymbol{F}}(T) \rightarrow \overline{\boldsymbol{F}}(T)$
  by setting

  \centerline{$ E^{\ast_{\iota}} := E^\ast\, \;\;\;\; \mbox{ for each } \;
   E \in \overline{\boldsymbol{F}}(T) \
  (\subseteq \overline{\boldsymbol{F}}(R) )\,.  $}

  \noindent Then we know \cite[Proposition 2.8]{FL1}

   \bf (a) \sl If $\iota$ is not the identity map, then  $\ast_{\iota}$
is a
      semistar, possibly non--star, operation on $T$, even if $\ast$ is
a star
      operation on $R$. \rm

     \rm Note that when $\ast$ is a star operation on $R$ and
     $(R:_{K}T) =(0)$, a fractional ideal $E$
     of $T$ is not necessarily a fractional ideal of $R$, hence
$\ast_{\iota}$
     is not defined as a star operation on $T$.

  \bf (b) \sl If $\ast$ is of finite type on $R$, then ${\ast_{\iota}} $
  is also of finite type on $T\,.$

  \bf (c) \sl When $T := R^\ast$, then ${\ast_{\iota}}$ defines a star
  operation on $T\,.$

  \bf (d) \sl If $\ast$ is e.a.b. [respectively, a.b.] on $R$ and if\ $T
:= R^\ast$, then ${\ast_{\iota}}$ is e.a.b.
  [respectively, a.b.] on $T$.

    \smallskip

    \rm Conversely, let $\star$ be a semistar operation on the overring
$T$ of
    $R$.  Define $\star^{\iota}: \overline{\boldsymbol{F}}(R)
\rightarrow
    \overline{\boldsymbol{F}}(R)$ by setting

  \centerline{$ E^{\star^{\iota}} := (ET)^\star\, \; \;\; \; \mbox{ for
each } \; E \in
  \overline{\boldsymbol{F}}(R)\,.
  $}

  \noindent Then we know \cite[Proposition 2.9, Corollary 2.10]{FL1}

  \bf (e) \sl ${\star^{\iota}} $ is a semistar operation on $R$.

  \bf (f) \sl If $\star:= d_{T}$, then ${(d_{T})}^{\iota}$ is a semistar
  operation of finite type on $R$, \rm  which is denoted also by $\star_{\{T\}}$
(i.e., it
  is the semistar operation on $R$ defined by $E^{\star_{\{T\}}} := ET$
  for each $E \in \oF(R)$).

  In particular, if $T=R$, then $\star_{\{R\}} = d_{R}$, and if
  $T=K$, then $\star_{\{K\}}= e_{R}$.  Note that if $R \subsetneq
  T$, then $\star_{\{T\}}$ is a proper semistar operation on $R$.

  \bf (g) \sl If $\star$ is e.a.b. [respectively, a.b.] on $T$, then
  ${\star^{\iota}}$ is e.a.b. [respectively, a.b.] on $R$.

  \bf (h) \sl For each semistar operation $\star$ on $T$, we have
  $({\star^{\iota}})_{\iota} = \star$.

  \bf (i) \sl For each semistar operation $\ast$ on $R$, we have
  $(\ast_{\iota})^{\iota} \geq \ast$ \rm (since
  $E^{({\ast_{\iota}})^{\iota}} = (ET)^{\ast_{\iota}} = (ET)^\ast
  \supseteq E^\ast$ for each $E \in \overline{\boldsymbol{F}}(R))$.

       \end{exxe}

       Other relevant classes of examples are recalled next.

       \begin{exxe} \label {ex:1.3}\rm Let $\Delta $ be a nonempty set
of prime
       ideals of an integral domain $R$ with quotient field $K$.  Set

  \centerline{ $E^{\star_{\Delta}}:=\cap \{ ER_P\; | \;\, P\in
  \Delta \}\, \;\;\;\;  \mbox{for each nonzero $R$--submodule $E$
  of $K$}\,.$}

  \noindent If $\Delta$ is the empty set, then we set $\star_\emptyset
:=
  e_{R}$.  The mapping $E\mapsto E^{\star_\Delta}$, for each $E\in
  \boldsymbol{\overline{F}}(R)$, defines a semistar operation on $R\,.$
\
  Moreover \cite [Lemma 4.1]{FH},

  \bf (a) \sl For each $E\in \boldsymbol{\overline{F}}(R)$ and for each
  $P\in \Delta\,$, $ER_P=E^{\star_\Delta}R_P$.

  \bf (b) \sl The semistar operation $\star_\Delta$ is \rm stable (with
  respect to the finite intersections), \sl i.e., for all $E, F \in
  \boldsymbol{\overline{F}}(R)\,$ we have $\, (E \cap
F)^{\star_{\Delta}}
  = E^{\star_{\Delta}} \cap F^{\star_{\Delta}}\,.  $

  \bf (c) \sl For each $P\in \Delta$, $P^{\star_\Delta}\cap R = P$.

  \bf (d) \sl For each nonzero integral ideal $I$ of $R$ such that
  $I^{\star_\Delta}\cap R\neq R$, there exists a prime ideal $P\in
\Delta$
  such that $I\subseteq P$.

\smallskip

  A \it semistar operation $\ast$ on $R$ \rm is called \it spectral \rm
if
  there exists a subset $\Delta $ of Spec($R$) such that $\ast =
  \star_{\Delta}\,;$ in this case, we say that $\ast$ is \it the
spectral
  semistar operation associated with \rm $\Delta \,.$

  We say that $\ast$ is a \it quasi--spectral semistar operation \rm
  (or that \it $\ast$ possesses enough primes\ \rm) if, for each
  nonzero integral ideal $I$ of $R$ such that $I^{\ast}\cap R\neq
  R$, there exists a prime ideal $P$ of $R$ such that $I\subseteq P$
  and $P^{\ast}\cap R = P\,.$ For instance, it is easy to see that
  \it if $\ast$ is a semistar operation of finite type, then $\ast$
  is quasi--spectral.  \rm

   From (c) and (d), we deduce that \sl each spectral semistar
  operation is quasi--spectral. \rm

  Given a semistar operation $\ast$ on $R$, assume that the set

  \centerline{$ \Pi^\ast:= \{ P \in \mbox{Spec}(R) \; | \;\; P\neq 0
\mbox{
  and } P^{\ast} \cap R \neq R \}$}

  \noindent is nonempty. Then the spectral semistar operation of $R$
  defined by $\, \ast_{sp}:= \star_{\Pi^{\ast}}\, $ is called \it
  the spectral semistar operation associated to $\ast\,.$ \rm \,
  Note that if $\ast$ is quasi--spectral such that $R^{\ast}\neq K$,
  then $ \Pi^{\ast}$ is nonempty and $ \ast_{sp} \leq \ast \, $
  \cite[Proposition 4.8 and Remark 4.9]{FH}.

  It is easy to see that \sl$\, \ast \,$ is spectral if and only if
  $\, \ast = \ast_{sp}\,.$ \rm

    For each semistar operation $\ast$ on $R$, we can consider

  \centerline{ $\tilde{\ast} := (\ast_{f})_{sp}\,.$}

  \noindent Then we know \cite[Proposition 3.6 (b), Proposition 4.23
(1)]{FH}

  \bf (e) \sl$\tilde{\ast}$ is a spectral semistar operation of finite
type
  on $R$, and \rm if $\calM(\ast_{f})$ denotes the set of all the
maximal
  elements in the set $\{I\mbox{ nonzero integral ideal of }$ $ R$ $ |
\;
  I^{\ast_{f}} \cap R \neq R \}\,,$ \ then \sl

  \centerline{$\tilde{\ast} = \star_{\calM{(\ast_{f})}}\,.$}

  \rm It is also known \cite[page 185]{FH} that for each $E \in
  \overline{\boldsymbol{F}}(R)$,

  \centerline{$E^{\tilde{\ast}} = \cup\{(E:_{K} F ) \mid \, F \in
\boldsymbol{f}
  (R)\,, \; F^{\ast}= R^\ast \}$\,.}

  \bf (f) \sl If $\ast$ is a star operation on $R$, then $\tilde{\ast}$
is
   a (spectral) star operation (of finite type) on $R$ and $\tilde{\ast}
\leq
   \ast$.  \rm

  If $\ast:= v_{R}$, using the notation introduced by Wang Fanggui and
  R.L. McCasland \cite{Fanggui/McCasland:1997}, we will denote by
$w_{R}$ (or simply by $w$) the star
  operation $\widetilde{v_{R}} = (t_{R})_{sp}$ (cf.  also
\cite{Hedstrom/Houston: 1980}
   and \cite{Anderson/Cook:2000}).

  \end{exxe}

  The construction of a spectral semistar operation associated to a set
of
  prime ideal can be generalized as follows.

  \begin{exxe} \label{ex:1.4} \rm Let $ \calR := \{ R_\lambda \; | \;\;
\lambda
  \in \Lambda \}$ be a nonempty family of overrings of $R\,$ and define
$\
  \star_{\calR}: \overline{\boldsymbol{F}}(R) \rightarrow
  \overline{\boldsymbol{F}}(R)$ by setting

  \centerline{$ E^{\star_{\calR}} := \cap \{ ER_\lambda \; | \;\;
\lambda
  \in \Lambda \}
  \, \; \mbox{ for each } \; E \in \overline{\boldsymbol{F}}(D) \,.  $}

  \noindent  Then we know \cite[ Lemma 2.4 (3), Example 2.5 (6),
Corollary
  3.8]{FL1}

  \bf (a) \sl The operation $ \star_{ {\calR}} $ is a semistar operation
  on $R$.  Moreover, if $ {\calR}= \{ R_P \; | \;\; P \in \Delta \},$
then
  $ \star_{\calR} = \star_{\Delta}\,.$

  \bf (b) \sl  $ E^{\star_{\calR} }R_\lambda = E R_\lambda \, \;$ for
each $E \in
  \overline{\boldsymbol{F}}(R)\,$ and for each $\lambda \in  \Lambda\,.$

  \bf (c) \sl If $ {\calR} = {\calW}$ is a family of valuation
  overrings of $R$, then $\star_{\calW}$ is an a.b. semistar
  operation on $D$.  \rm

  \smallskip

  We say that \it two semistar operations on \rm $D$, $\star_{1}$ and
$\star_{2}$ are \it
  equivalent \rm if $(\star_{1})_{f} = (\star_{2})_{f}\,.$
  Then we know (\cite[Proposition 3.4]{FL2} and \cite[Theorem 32.12]{G})

  \bf (d) \sl Each e.a.b. semistar [respectively, star] operation on $R$
  is equivalent to a  semistar [respectively, star] operation of the
type
  $\star_{\calW}$ for some family $\calW$ of valuation overrings of $R$
  [respectively, for some family $\calW$ of valuation overrings of $R$
  such that $R = \cap \{ W \mid W \in \calW\}$].  \rm

  If \ ${\calW}$ \ is the family of \it all \rm the valuation overrings
of
  $R$, then $\star_{\calW}$ is called \it the $b_{R}$--semistar
operation
  \rm (or simply  \it the $b$--semistar operation on $R$\rm ).\
  Moreover, if $R$ is integrally closed, then\ $R^{b_{R}} = R\,$
  \cite[Theorem 19.8]{G}, and thus the operation $b $ defines a star
  operation on $R\,,$ which is called \it the $b$-star operation \rm \cite[p.
  398]{G}.  \rm

  \end{exxe}

  \begin{exxe} \label{ex:1.5} \rm If $\{\ast_{\lambda}\mid \lambda \in
  \Lambda\}$ is a family of semistar [respectively, star] o\-pe\-rations
  on $R$, then \ $\wedge_{\lambda} \{\ast_{\lambda}\mid \lambda \in
  \Lambda\}$ \ (denoted simply by \ $\wedge \ast_{\lambda}$\ ), defined
by

       \centerline{
       $ E^{\wedge \ast_{\lambda}}:= \cap \{ E^{\ast_{\lambda}}\mid
\lambda \in
       A\}\,$ \;\;\;\ {\mbox{for each $E \in
\overline{\boldsymbol{F}}(R)$ \
       [respectively, $E \in \boldsymbol{F}(R)$]\,,} } }

    \noindent is a semistar [respectively, star] operation on $R$.  This
type
    of semistar o\-pe\-ra\-tion generalizes the semistar [respectively,
    star] operation of type $\star_{\calR}$ (where $ \calR := \{
R_\lambda
    \mid \lambda \in \Lambda \}$ is a nonempty family of overrings of
$R$;
    Example \ref{ex:1.4}), since

    \centerline{$\star_{\calR} = \wedge\
    \star_{\{R_{\lambda}\}}$}

    \noindent where $\star_{\{R_{\lambda}\}}$ is the semistar operation
on $R$
    considered in Example \ref{ex:1.2} (f).

    Note the following observations:

    \bf (a) \sl If at least one of the semistar operations in the family
    $\{\ast_{\lambda}\mid \lambda \in \Lambda\}$ is a star operation on
$R$,
    then\ $\wedge \ast_{\lambda}$\ is still a star operation on $R$.\rm

    \bf (b) \sl  Let $\iota:R \hookrightarrow T$ be an embedding of
integral
    domains with the same field of quotients $K$ and let
$\{\ast_{\lambda}\mid
    \lambda \in \Lambda\}$ be a family of semistar operations on $R$.
Then

   \vskip - 6pt  \centerline{$ (\wedge \ast_{\lambda})_{\iota} = \wedge
    (\ast_{\lambda})_{\iota}$\,.}

   \bf (c) \sl  Let $\iota:R \hookrightarrow T$ be an embedding of
integral
    domains with the same field of quotients $K$ and let
$\{\star_{\lambda}\mid
    \lambda \in \Lambda\}$ be a family of semistar operations on $T$,
then

    \vskip - 6pt  \centerline{$ (\wedge \star_{\lambda})^{\iota} =
\wedge
    (\star_{\lambda})^{\iota}$\,.} \rm

    \end{exxe}

   \smallskip

  \section{Star operations and pullbacks}

  For the duration of this paper we will mainly consider the following
  situations:

  \smallskip
  \setlength{\hangindent}{23pt} \noindent \bf (\th) \sl  $\,T\,$
  represents an integral domain,\ $M$\ an ideal of \ $T$,\ $k$\ the
  factor ring\ $T/M$,\ $D$\ an integral domain subring of\ $k$\  and
$\varphi: T \rightarrow T/M =: k$ the canonical
  projection.  Set $R: = \varphi^{-1}(D)=: T\times_{k }D$ \rm the
pullback
  of $D$ inside $T$ with respect to $\varphi$\sl , hence $R$ is an
  integral domain (subring of $T$).  Let $K$ denote the field of quotients of $R$.

  \smallskip
  \setlength{\hangindent}{23pt} \noindent  \bf (\th$^{+}$) \sl Let
  $L$ be the field of quotients of $D$.
  In the situation \bf (\th)\sl, we assume, moreover,
  that $L \subseteq k$, and denote by $S:=
  \varphi^{-1}(L)=: T\times_{k }L$ \rm the pullback of $L$ inside $T$
with respect to $\varphi$\sl.  Then $S$ is an integral domain with field of
  quotients equal to $K$.  In this situation, $M$, which is a prime
ideal in $R$, is a maximal ideal in $S$.  Moreover, if $M\neq (0)$ and
  $D\subsetneq k$, then $M$ is a divisorial ideal of $R$, actually,
  $M=(R:T)$.

  \vskip.1in  \rm

  Let $\star_{D}$ [respectively, $\star_{T}$] be a star operation on
  the integral domain $D$ [respectively, $T$].  Our first goal is to
  define in a natural way a star operation on $R$, which we will
  denote by {\scriptsize $\Diamond$}, associated to the given star
  operations on $D$ and $T$.  More precisely, if we denote by\ \bf
  Star\rm$(A)$ \ the set of all the star operations on an integral
  domain $A$, then we want to define a map

  \centerline{\, $\Phi:$ \bf Star\rm$(D)\ \times $ \bf Star\rm$(T)
\rightarrow$ \bf
  Star\rm$(R)$, \,\; $(\star_{D}, \star_{T}) \mapsto
\mbox{\scriptsize{$\Diamond$}}$ \,.}

  For each nonzero fractional ideal $I$ of $R$, set
  $$
   I^{\mbox{\tiny{$\Diamond$}}} := \cap
  \left\{
  x^{-1}\varphi^{-1}\left(\left(\frac{xI+M}{M}\right)^{\star_{D}}\right)
  \mid \, x \in I^{-1}\,, \; x \neq 0 \right\} \, \cap \,
  (IT)^{\star_{T}}\,,
  $$
  where if $\frac{xI+M}{M}$ is the zero ideal of $D$ (i.e., if $xI+M
  \subseteq M$), then we set
  $\varphi^{-1}\left(\left(\frac{xI+M}{M}\right)^{\star_{D}}\right) =
M$.

  \vskip .1in

  \begin{prro}\label{pr:2.1} Keeping the notation and hypotheses
  introduced in \rm \bf (\th)\it , then $\mbox{\scriptsize{$\Diamond$}}$
  defines a star operation on the integral domain $R$\ ($=T\times_{k
}D$).
  \end{prro}
   \rm

  \noindent {\bf Proof}. {\bf Claim 1}.  \sl For each nonzero
  fractional ideal $I$ of $R$, $I
   \subseteq I^{\mbox{\tiny{$\Diamond$}}}$.\rm

  We have
  $$
  \begin{array}{rl}
      I^{\mbox{\tiny{$\Diamond$}}} &\supseteq \cap \left\{
      x^{-1}\varphi^{-1}\left(\frac{xI+M}{M}\right) \mid \, x \in
I^{-1}\,, \;
      x \neq 0 \right\} \, \cap \, IT\,\\
  & = \cap \left\{ x^{-1}\left(xI+M\right) \mid \, x \in I^{-1}\,,
  \;  x \neq 0  \right\} \, \cap \, IT\ \\ & = \cap \left\{
  I+x^{-1}M \mid \, x \in I^{-1}\,, \;  x \neq 0  \right\} \, \cap
  \, IT\, \supseteq I\,.
  \end{array}$$

  \smallskip

   {\bf Claim 2}.  \sl For each nonzero element $z$ of $K$,
   $(zR)^{\mbox{\tiny{$\Diamond$}}} = zR$ (in parti\-cu\-lar,
$R^{\mbox{\tiny{$\Diamond$}}} = R$).\rm

  We have
  $$
  \begin{array}{rl}
      (zR)^{\mbox{\tiny{$\Diamond$}}} &= \cap \left\{

x^{-1}\varphi^{-1}\left(\left(\frac{xzR+M}{M}\right)^{\star_{D}}\right)
      \mid \, x \in z^{-1}R\,, \; x \neq 0 \right\} \, \cap \,
      (zT)^{\star_{T}}  \\
  & \subseteq

z\left(\varphi^{-1}\left(\left(\frac{R+M}{M}\right)^{\star_{D}}\right)\right)

  \, \cap \, zT\, \\ & =
  z\left(\varphi^{-1}\left(\frac{R}{M}\right)\right) \, \cap \, zT\,
  = zR \cap zT = zR \,.
  \end{array}$$
  Therefore, by Claim 1, we deduce that $(zR)^{\mbox{\tiny{$\Diamond$}}}
= zR$.

  \smallskip

  {\bf Claim 3}.  \sl For each nonzero element $z$ of $K$ and for
  each nonzero fractional ideal $I$ of $R$,
  $(zI)^{\mbox{\tiny{$\Diamond$}}} =
  zI^{\mbox{\tiny{$\Diamond$}}}$.\rm

  Note that given $0 \neq z \in K$, for each nonzero $x  \in I^{-1}$
   there exists a unique $y \in (zI)^{-1} $ such that $x =yz$.  Therefore, we have
  $$
  \begin{array}{rl}
      I^{\mbox{\tiny{$\Diamond$}}} &= \cap \left\{

x^{-1}\varphi^{-1}\left(\left(\frac{xI+M}{M}\right)^{\star_{D}}\right)
      \mid \, x \in I^{-1}\,, \; x \neq 0 \right\} \, \cap \,
      (IT)^{\star_{T}}\, \\
  &= \cap \left\{

(yz)^{-1}\varphi^{-1}\left(\left(\frac{yzI+M}{M}\right)^{\star_{D}}\right)

  \mid \, yz \in  I^{-1}\,, \;  y \neq 0 \right\} \, \cap \,
  (IT)^{\star_{T}} \\ &= \cap \left\{

z^{-1}y^{-1}\varphi^{-1}\left(\left(\frac{yzI+M}{M}\right)^{\star_{D}}\right)

  \mid \, y \in (zI)^{-1} \,,\; y \neq 0 \right\} \, \cap \,
  (IT)^{\star_{T}} \\ &= z^{-1}\left( \cap \left\{

y^{-1}\varphi^{-1}\left(\left(\frac{yzI+M}{M}\right)^{\star_{D}}\right)
  \mid \, y \in (zI)^{-1} \,,\; y \neq 0 \right\}\, \cap \,
  (zIT)^{\star_{T}}\right) \\ &=
  z^{-1}(zI)^{\mbox{\tiny{$\Diamond$}}}\,.
  \end{array}
  $$
  Thus  we immediately conclude that $(zI)^{\mbox{\tiny{$\Diamond$}}} =
  zI^{\mbox{\tiny{$\Diamond$}}}$.

  \smallskip

  {\bf Claim 4}.  \sl For each pair of nonzero fractional ideals
  $I\subseteq J$ of $R$, $I^{\mbox{\tiny{$\Diamond$}}} \subseteq
  J^{\mbox{\tiny{$\Diamond$}}}$.\rm

  Since $J^{-1} \subseteq I^{-1}$, we have
  $$
  \begin{array}{rl}
  J^{\mbox{\tiny{$\Diamond$}}} &= \cap \left\{
  x^{-1}\varphi^{-1}\left(\left(\frac{xJ+M}{M}\right)^{\star_{D}}\right)

  \mid \, x \in J^{-1}\,, \; x \neq 0 \right\} \, \cap \,
  (JT)^{\star_{T}} \\ &\supseteq  \cap \left\{
  x^{-1}\varphi^{-1}\left(\left(\frac{xI+M}{M}\right)^{\star_{D}}\right)

  \mid \,  x \in J^{-1}\,, \;  x \neq 0 \right\} \, \cap \,
  (IT)^{\star_{T}}  \\ &\supseteq  \cap \left\{
  x^{-1}\varphi^{-1}\left(\left(\frac{xI+M}{M}\right)^{\star_{D}}\right)

  \mid \, x \in I^{-1}\,, \;  x \neq 0  \right\} \, \cap \,
  (IT)^{\star_{T}} = I^{\mbox{\tiny{$\Diamond$}}} \,.
  \end{array}
  $$

  \smallskip

  {\bf Claim 5}.  \sl For each nonzero fractional ideal $I$ of $R$,
   $I \subseteq I^{{\mbox{\tiny{$\Diamond$}}}} \subseteq I^v$,
  and hence $(I^{{\mbox{\tiny{$\Diamond$}}}})^{-1} = I^{-1}$.  \rm

  Since $ I^v = \cap \{zR \mid\, I \subseteq zR\,,\; z \in K\}$,
  by Claim 2, we deduce that

  \centerline{$ I \subseteq zR \,\; \Rightarrow \,\;
  I^{{\mbox{\tiny{$\Diamond$}}}} \subseteq
  (zR)^{{\mbox{\tiny{$\Diamond$}}}} = zR\,,
  $}

  \noindent hence $I^{\mbox{\tiny{$\Diamond$}}} \subseteq I^v$.

  \smallskip

  {\bf Claim 6}.  \sl For each nonzero fractional ideal $I$ of $R$,
  $(I^{\mbox{\tiny{$\Diamond$}}})^{\mbox{\tiny{$\Diamond$}}} =
  I^{\mbox{\tiny{$\Diamond$}}}$.\rm

  Since $(I^{\mbox{\tiny{$\Diamond$}}})^{-1} = I^{-1}$ for each
  nonzero ideal $I$ of $R$, we have $$
  \begin{array}{rl}
  (I^{\mbox{\tiny{$\Diamond$}}})^{\mbox{\tiny{$\Diamond$}}} &= \cap
  \left\{

x^{-1}\varphi^{-1}\left(\left(\frac{xI^{\mbox{\tiny{$\Diamond$}}}+M}{M}\right)^{\star_{D}}\right)

  \mid \, x \in (I^{\mbox{\tiny{$\Diamond$}}})^{-1}\,, \; x \neq 0
  \right\} \, \cap \,
  (I^{{\mbox{\tiny{$\Diamond$}}}}T)^{\star_{T}}\, \\ &= \cap \left\{

x^{-1}\varphi^{-1}\left(\left(\frac{xI^{\mbox{\tiny{$\Diamond$}}}+M}{M}\right)^{\star_{D}}\right)

  \mid \, x \in I^{-1}\,, \; x \neq 0 \right\} \, \cap \,
  (I^{{\mbox{\tiny{$\Diamond$}}}}T)^{\star_{T}}\, .
  \end{array}
  $$

Note that for $0 \neq x\in I^{-1}$ with $xI\subseteq M$, we
have \smallskip

\hskip 0.6cm $\bullet \hskip .3cm
xI^{\mbox{\tiny{$\Diamond$}}}\subseteq M$\; (\mbox{and so }\,
$x^{-1}\varphi^{-1}\left(\left(\frac{xI^{\mbox{\tiny{$\Diamond$}}}+M}{M}\right)^{\star_{D}}\right)=
x^{-1}\varphi^{-1}\left(\left(\frac{xI+M}{M}\right)^{\star_{D}}\right)$),
since
$I^{\mbox{\tiny{$\Diamond$}}} \subseteq
x^{-1}\varphi^{-1}\left(\left(\frac{xI+M}{M}\right)^{\star_{D}}\right)=x^{-1}M$.

   Now for $0\neq x\in I^{-1}$ with $xI\not\subseteq M$, we
   have \smallskip

  \hskip 0.6cm $\bullet \hskip .3cm
  \left(\frac{xI^{\mbox{\tiny{$\Diamond$}}}+M}{M}\right)^{\star_{D}} \,
\subseteq
  \left(\frac{xI+M}{M}\right)^{\star_{D}}\,\; \mbox{ (and so } \,
  \left(\frac{xI^{\mbox{\tiny{$\Diamond$}}} + M}{M}\right)^{\star_{D}} =
  \left(\frac{xI+M}{M}\right)^{\star_{D}} )\,,$
  \smallskip

\noindent since
  $$
  \begin{array}{rl}
   &I^{\mbox{\tiny{$\Diamond$}}} \subseteq

x^{-1}\varphi^{-1}\left(\left(\frac{xI+M}{M}\right)^{\star_{D}}\right)
   \;\;  \Rightarrow \; \; xI^{\mbox{\tiny{$\Diamond$}}} \subseteq
      \varphi^{-1}\left(\left(\frac{xI+M}{M}\right)^{\star_{D}}\right)
\;\;
       \\
   &\;\; \Rightarrow \; \frac{xI^{\mbox{\tiny{$\Diamond$}}}+M}{M} =
   \varphi\left(xI^{\mbox{\tiny{$\Diamond$}}}\right) \, \subseteq
\varphi\left(
   \varphi^{-1}\left(\left(\frac{xI+M}{M}\right)^{\star_{D}}\right)
\right)
   = \left(\frac{xI+M}{M}\right)^{\star_{D}} \, \;\;  \\
       &\;\; \Rightarrow \;
\left(\frac{xI^{\mbox{\tiny{$\Diamond$}}}+M}{M}\right)^{\star_{D}} \,
       \subseteq \left(\frac{xI+M}{M}\right)^{\star_{D}}\,.
       \end{array}
  $$

  Lastly,

  \hskip 0.6cm $\bullet \hskip .3cm
(I^{{\mbox{\tiny{$\Diamond$}}}}T)^{\star_{T}}\subseteq
  (IT)^{\star_{T}}\, (\mbox{and so }\,
(I^{{\mbox{\tiny{$\Diamond$}}}}T)^{\star_{T}}=
  (IT)^{\star_{T}})\,,$

  \noindent since

  \centerline{$ I^{\mbox{\tiny{$\Diamond$}}} \subseteq (IT)^{\star_{T}}
  \Rightarrow I^{{\mbox{\tiny{$\Diamond$}}}}T \subseteq (IT)^{\star_{T}}
  \Rightarrow (I^{{\mbox{\tiny{$\Diamond$}}}}T)^{\star_{T}} \subseteq
  (IT)^{\star_{T}} \,.
      $}

  Therefore, we can easily  conclude
  $$
  \begin{array}{rl}
      (I^{\mbox{\tiny{$\Diamond$}}})^{\mbox{\tiny{$\Diamond$}}} &= \cap
\left\{

x^{-1}\varphi^{-1}\left(\left(\frac{xI^{\mbox{\tiny{$\Diamond$}}}+M}{M}\right)^{\star_{D}}\right)

  \mid \, x \in I^{-1}\,, \; x \neq 0 \right\} \, \cap \,
  (I^{{\mbox{\tiny{$\Diamond$}}}}T)^{\star_{T}}  \\ &= \cap \left\{
  x^{-1}\varphi^{-1}\left(\left(\frac{xI+M}{M}\right)^{\star_{D}}\right)

  \mid \, x \in I^{-1}\,, \; x \neq 0 \right\} \, \cap \,
  (IT)^{\star_{T}} = I^{\mbox{\tiny{$\Diamond$}}} \,.
       \end{array}
  $$

  The previous argument shows that $\mbox{\scriptsize{$\Diamond$}}$ is a
(well defined) star
  operation on the integral domain $R$.  \hfill $\Box$

  \vspace{.1in}

  \begin{reem} \label{rk:2.2} \rm
      \bf (a) \rm Note that in the
  proof of Proposition \ref{pr:2.1} $M$ is possibly a nonmaximal ideal
of $T$ (and
  $R$), even though we assume that $M \,(= M \cap R)$ is a prime ideal
of
  $R$.

  \bf (b) \rm In the pullback setting {\bf (\th)}, for each nonzero
ideal
  $I$ of $R$ with $I \subseteq M$, $I^{\diamond}\subseteq M$, because
  $I^{\diamond}\subseteq
  x^{-1}{\varphi}^{-1}((\frac{xI+M}{M})^{{\star}_D})$ for each $x\in
  I^{-1}\setminus (0)$ and $1\in R\subseteq I^{-1}$, thus
  $I^{\diamond}\subseteq {\varphi}^{-1}((\frac{I+M}{M})^{{\star}_D})=M$.
  In particular, if $M\neq (0)$, then $M= M^{\diamond}$.

  \bf (c) \rm If $M$ is the ideal $(0)$, then $T=k$ and $R =D$.  In this
  extreme situation, we have
  $\mbox{\scriptsize{$\Diamond$}} = \star_{D} \wedge
({\star}_T)^{\iota}$,
  where $\iota: R=D \hookrightarrow T=k$ is the canonical inclusion.

  \noindent Note that it can happen that
$\mbox{\scriptsize{$\Diamond$}} =
  \star_{D} \wedge ({\star}_T)^{\iota} \lneq \star_{D}$.  For instance,
  let $R$ be a Krull domain of dimension $\geq 2$, $P$ a prime ideal of
  $R$ with ht$(P) \geq 2$, $T:=R_{P}$ and $M := (0)$ (hence, $R=D$ and
  $T=k$).  Set $\star_{D}:=v_{D}$ and $\star_{T} := d_{T}$. Then
$P^{\star_{D}}
  =P^{v_{D}} =P^{v_{R}} =R$, but $P^{\diamond} = P^{v_{D}} \cap
  (PT)^{\star_{T} }= R \cap PR_{P} = P$.

  \bf (d) \rm Let $M \neq (0)$. If $D =  L=k$ (in particular, $M$ must
be a
  nonzero maximal ideal of $T$, and necessarily, $\star_{D}$ is the
(unique)
  star operation $d_{D}$ of $D=L=k$), then $R = T$.  In this extreme
  situation, we have that $\mbox{\scriptsize{$\Diamond$}}$ and
$\star_{T}$
  are two star operations on $T$\ (with $\mbox{\scriptsize{$\Diamond$}}
  \leq \star_{T}$)\ that are possibly different.  For instance, if $T$
is an
  integral domain with a nonzero nondivisorial maximal ideal $M$ (e.g.
  $T:=k[X, Y]$, $M:= (X,Y)$) and if $\star_{T} := v_{T}$, then
  $M^{\diamond}=M$  by (b), but $M^{\star_{T}}= M^{v_{T}} = T$.

      If $D=k$, but $D\subsetneq L$, then it is not difficult to see
that
      $\mbox{\scriptsize{$\Diamond$}} = \star_{T}$ if and only if, for
each
      nonzero ideal $I$ of $R=T$ with $I \not\subseteq M$,
$\frac{I^{\star_{T}}
      +M}{M} \subseteq \left(\frac{I+M}{M}\right)^{\star_{D}}$.

  \end{reem}

  \vspace{.02in}

  Our next example will explicitly show the behaviour of the star
operation
  $\mbox{\scriptsize{$\Diamond$}}$ in some special cases of the pullback
  construction \rm \bf (\th)\rm .

  \vspace{.02in}

  \begin{exxe} \label{ex:2.3} \sl With the notation and hypotheses
  introduced in \rm \bf (\th)\sl, assume, moreover, that $T$ is local
with nonzero
  maximal ideal $M$ and $D=L$ is a proper subfield of $k$. \rm In this
  special case of the situation \rm \bf (\th$^{+}$)\rm ,
  $\star_{D}=d_{D}=e_{D}$ is the unique star operation on $D$.  \sl Let
$I$ be
  a nonzero fractional ideal of $R$.

  \balf \bf \item \sl If $II^{-1} = R$, then
$I^{\mbox{\tiny{$\Diamond$}}}
  = I = I^{v_{R}}$.  \rm

  \bf \item \sl If $II^{-1} \subsetneq R$, then
  $I^{\mbox{\tiny{$\Diamond$}}} = I^{v_{R}}\cap (IT)^{\star_{T}}$.
  Moreover, if $(IT)^{\star_{T}} = x^{-1}T$ for some nonzero $x \in
  I^{-1}$, then $I^{\mbox{\tiny{$\Diamond$}}} = I^{v_{R}}\ (\subsetneq
  (IT)^{\star_{T}})$.  If $(IT)^{\star_{T}} \neq x^{-1}T$ for all $x \in
  I^{-1}$, then $I^{\mbox{\tiny{$\Diamond$}}} = (IT)^{\star_{T}}$.  \rm

  \bf \item \sl If $[k:L] > 2$ and if $T \subsetneq M^{-1} =(R:_{K} M)$,
  then $d_{R} \neq \mbox{\scriptsize{$\Diamond$}}\neq v_{R}$ for all the
star operations
  $\star_{T}$ on $T$.  \rm

  \bf \item \sl Let $[k:L] = 2$.  If $T$ is (local but) not a valuation
  domain, then $d_{R} \neq \mbox{\scriptsize{$\Diamond$}}$ for all the
star operations
  $\star_{T}$ on $T$.  If $T = (R:_{K}M)$ and if $\star_{T}=v_{T}$, then
  $\mbox{\scriptsize{$\Diamond$}} = v_{R}$.  \rm

  \ealf

  \bf (a) \rm is obvious, because $I$ is invertible, hence $I$ is
  divisorial (in fact, $I$ is principal, since $R$ is also local) and so
$I =
  I^{\mbox{\tiny{$\Diamond$}}} = I^{v_{R}} \ (\subseteq
(IT)^{\star_{T}})$.

  \bf (b) \rm Note that for each nonzero ideal $I$ of $R$ with the
  property that $II^{-1} \subsetneq R$, we have necessarily that
$II^{-1}
  \subseteq M$.  Moreover, for each nonzero $x \in I^{-1}$, from $xI
  \subseteq M$, we deduce that $I \subseteq x^{-1}M$ and so we have that
  $I^{v_{R}}= \cap \{x^{-1}M\mid \, x \in I^{-1}\,,\; x \neq 0 \}$.
  Therefore,
   $$
  \begin{array}{rl}
      I^{\mbox{\tiny{$\Diamond$}}} &= \cap \left\{

x^{-1}\varphi^{-1}\left(\left(\frac{xI+M}{M}\right)^{\star_{D}}\right)
      \mid \, x \in I^{-1}\,, \; x \neq 0 \right\} \, \cap \,
      (IT)^{\star_{T}} \\
  &= \cap \left\{
  x^{-1}\varphi^{-1}\left(\left(\frac{M}{M}\right)^{\star_{D}}\right)
  \mid \, x \in I^{-1}\,, \;  x \neq 0  \right\} \, \cap \,
  (IT)^{\star_{T}} \\ &= \cap \left\{ x^{-1}M \mid \, x \in
  I^{-1}\,, \; x \neq 0 \right\} \, \cap \, (IT)^{\star_{T}} =
  I^{v_{R}} \, \cap \, (IT)^{\star_{T}}\,.
  \end{array}
      $$

      In order to prove the second part of \bf (b)\rm , note that in
this
      case, for each $0 \neq x \in I^{-1}$, we have
      $$
  \begin{array}{rl}
      I \subsetneq x^{-1}R \, &\Rightarrow \, I \subseteq x^{-1}M \,
\Rightarrow \,
      IT \subseteq x^{-1}MT=x^{-1}M \,  \\
      &\Rightarrow \, (IT)^{\star_{T}} \subseteq (x^{-1}M)^{\star_{T}} =

      x^{-1}M^{\star_{T}}\subseteq x^{-1}T\,.
      \end{array}
      $$

      Therefore, if  $(IT)^{\star_{T}} = x^{-1}T$ for some nonzero
  $x \in I^{-1}$, then $I^{v_{R}} \subseteq x^{-1}R \subseteq
  x^{-1}T = (IT)^{\star_{T}}$.  Thus, in this case,  $
  I^{\mbox{\tiny{$\Diamond$}}} = I^{v_{R}} $. Assume that
  $(IT)^{\star_{T}} \subsetneq x^{-1}T$ for all $x \in I^{-1}$.
  Then $(IT)^{\star_{T}} \subseteq x^{-1}M$ and thus

  \centerline{$ (IT)^{\star_{T}} \subseteq \cap \left\{
  x^{-1}M \mid \, x \in I^{-1}\,, \; x \neq 0 \right\} = I^{v_{R}}\,,$}

  \noindent hence $I^{\mbox{\tiny{$\Diamond$}}} = (IT)^{\star_{T}}$.

  \bf (c) \rm Let $0 \neq a\in M$, and let $z\in T\setminus R$.  Set
  $I:= (a, az)R$. Then obviously $IT =aT$ (since $z$ is invertible in
$T$), thus
  $(IT)^{\star_{T}}=aT=IT$.

  Note that, in this case, $(IT)^{\star_{T}}=aT \subsetneq x^{-1}T$
  for all $x \in I^{-1}$. As a matter of fact, if $aT = x^{-1}T$ for
  some $x \in I^{-1}$, then $ax =u$ is a unit in $T$ and $ax \in R$
  (because $a\in I$ and $x \in I^{-1}$).  Hence, $ax$ is a unit in $R$.
  Now we reach a contradiction, since we deduce that $I \subseteq
x^{-1}R
  = aR \subseteq I$, i.e., $I = aR$.

  By (b), we have that $I^{\mbox{\tiny{$\Diamond$}}} =
  (IT)^{\star_{T}}=aT =IT \supsetneq I$, hence
  $\mbox{\scriptsize{$\Diamond$}} \neq d_{R}$.

  Assume also that $T\subsetneq M^{-1}$.  Since
  $I^{\mbox{\tiny{$\Diamond$}}} = (IT)^{\star_{T}}=aT$,  $I^{v_{R}}
  = (I^{\mbox{\tiny{$\Diamond$}}})^{v_{R}} = (aT)^{v_{R}}=
  a(R:_{K}(R:_{K}T)) = a(R:_{K}M)\supsetneq aT =
  (IT)^{\star_{T}}=I^{\mbox{\tiny{$\Diamond$}}}$.  Therefore
  $\mbox{\scriptsize{$\Diamond$}} \neq v_{R}$.

   \bf (d) \rm In the present situation, we can find $a, b \in M$ such
that
   $aT \not\subseteq bT$ and $bT \not\subseteq aT$.  Set $I :=(a,b)R$.

   It is easy to see that $I$ is not a principal ideal of $R$. (If
   $I =(a,b)R =cR$, then $a =cr_{1}, b =cr_{2}, c =as_{1}+bs_{2}$ and so
   $1=r_{1}s_{1}+r_{2}s_{2}$\ for some $r_{1}, s_{1}, r_{2}, s_{2} \in
R$;
   hence either $r_{1}s_{1}$ or $r_{2}s_{2}= 1-r_{1}s_{1}$ is a unit in
the
   local ring $R$.  For instance, if $r_{1}s_{1}$ is a unit in $R$, then
   $r_{1}$ is also a unit in $R$ and so $cR = aR$.  Thus $bR \subseteq
aR$,
   contradicting the choice of $a$ and $b$.)

   Note that $I$ is not a divisorial ideal of $R$.
    As a matter of fact, if $I = I^{v_{R}}$, then $I$ should be also an
   ideal of $T$ (i.e., $I = IT$) by \cite[Corollary
2.10]{Gabelli/Houston}.
   On the other hand, if $ z\in T\setminus R$, then $az\in IT =I
=(a,b)R$
   and so $az =ar_{1}+br_{2}$, i.e., $a(z-r_{1})=br_{2}$ for some $r_{1},
   r_{2}\in R$.  If $z-r_{1}\in M$, then $z\in r_{1}+M \subseteq R$,
which
   contradicts the choice of $z$.  If $z-r_{1}\in T\setminus M,$ then $a
   =br_{2}(z-r_{1})^{-1}\in bT$, which contradicts the choice of $a$ and
   $b$.  Hence, $I\neq IT$ and so $I \neq I^{v_{R}}$.

   If $(IT)^{\star_{T}} = x^{-1}T$ for some nonzero $x \in
   I^{-1}$, then (by (b)) $I^{\mbox{\tiny{$\Diamond$}}} = I^{v_{R}}\neq
I$, and so $d_{R}\neq
   \mbox{\scriptsize{$\Diamond$}} $.
   Assume that $(IT)^{\star_{T}} \neq x^{-1}T$ for all $x \in
    I^{-1}$, then (by (b)) $I^{\mbox{\tiny{$\Diamond$}}} =
(IT)^{\star_{T}}\supseteq IT
    \supsetneq I$, and so $d_{R}\neq \mbox{\scriptsize{$\Diamond$}} $.

   Finally, suppose that $T =(R:_{K}M)$ and that $\star_{T} = v_{T}$.
  Let $J$ be a nonzero fractional ideal of $R$. If $J$ is
  divisorial, then obviously $J^{\diamond}=J=J^{v_R}$. Assume that
  $J$ is not divisorial, then $JJ^{-1}\subsetneq R$.
   If $(JT)^{\star_{T}} = x^{-1}T$ for some nonzero $x \in
   J^{-1}$, then (by (b)) $J^{\mbox{\tiny{$\Diamond$}}} = J^{v_{R}}$.
   If $(JT)^{v_{T}} \neq x^{-1}T$ for all $x \in
    J^{-1}$,  then (by (b)) $J^{\mbox{\tiny{$\Diamond$}}} =
(JT)^{v_{T}}$.  Since $T =(R:_{K}M)=
    (M:_{K}M)$,  every  divisorial ideal of $T$ is  divisorial as an
    ideal of $R$ by \cite[Corollary 2.9]{Gabelli/Houston}.  Therefore

    \centerline{$J^{v_{R}} = (J^{\mbox{\tiny{$\Diamond$}}} )^{v_{R}} =
    ((JT)^{v_{T}})^{v_{R}}= (JT)^{v_{T}}
    =J^{\mbox{\tiny{$\Diamond$}}}\,,$}

    \noindent hence we conclude that $\mbox{\scriptsize{$\Diamond$}} =
{v_{R}}$.
      \end{exxe}

  \vspace{4pt} The previous construction of the star operation
  $\mbox{\scriptsize{$\Diamond$}}$ on the integral domain $R$ arising
from a
  pullback diagram gives the idea for ``lifting a star operation'' with
  respect to a surjective ring homomorphim between two integral domains.

  \begin{coor} \label{cr:2.4} Let $R$ be an integral domain with field
of
  quotients $K$, $M$ a prime ideal of $R$. Let $D$ be the
  factor ring $R/M$ and let $\varphi: R \rightarrow D$ be the canonical
  projection.  Assume that $\star$ is a star operation on $D$.  For each
  nonzero fractional ideal $I$ of $R$, set $$\begin{array} {rl}
  I^{\star^{\varphi}} &:= \cap \left\{
  x^{-1}\varphi^{-1}\left(\left(\frac{xI+M}{M}\right)^{\star}\right)
\mid
  \, x \in I^{-1}\,, \; x \neq 0 \right\}\\ &= \cap \left\{
  x\varphi^{-1}\left(\left(\frac{x^{-1}I+M}{M}\right)^{\star}\right)
\mid
  \, x \in K\,, \; I \subseteq xR \right\}\,,
  \end{array}
  $$
  where, as before, if\, $\frac{zI+M}{M}$\, is the zero ideal of $D$,
  then we set
  $\varphi^{-1}\left(\left(\frac{zI+M}{M}\right)^{\star}\right)$ = $M$.
  Then ${\star^{\varphi}}$ is a star operation on $R$.

      \end{coor}
       \rm

  \noindent {\bf Proof}.  \it Mutatis mutandis \rm the arguments used in
  the proof of Proposition~\ref{pr:2.1} show that ${\star^{\varphi}}$ is
a star
  operation on $R$.  \hfill $\Box$

     \vspace{10pt}

    Using the notation introduced in Section 1, in particular,
    in Example \ref{ex:1.2}, we immediately have the following:


  \begin{coor}\label{cr:2.5} With the notation and hypotheses introduced
in \rm
  \bf (\th)\it\ and Proposition \ref{pr:2.1}, if we use the definition
given
  in Corollary \ref{cr:2.4}, we have

  \centerline{$ \mbox{\scriptsize{$\Diamond$}} = (\star_{D})^\varphi
\wedge
  (\star_{T})^\iota\,.
  $}
  \end{coor}
  \vskip -20 pt \hfill $\Box$

  \vspace{12pt}

 We next examine the problem of ``projecting a star operation'' with
  respect to a surjective homomorphism of integral domains.

  \begin{prro}\label{pr:2.6} Let $R$, $K$, $M$, $D$, $\varphi$ be as in
Corollary
  \ref{cr:2.4} and let $L$ be the field of quotients of $D$.  Let $\ast$
  be a given star operation on the integral domain $R$.  For each
nonzero
  fractional ideal $F$ of $D$, set

  \centerline{$ F^{\ast_{\varphi}}:= \cap
\left\{y\varphi\left(\left(\varphi^{-1}\left(y^{-1}F\right)\right)^{\ast}\right)
  \mid \, y \in L\,,\;  F \subseteq yD \right\} \,.$}

  \noindent Then $\ast_{\varphi}$ is a star operation on $D$.
  \end{prro}
   \rm

  \noindent {\bf Proof}.  The following claim is a straightforward
  consequence of the definition.

  {\bf Claim 1}.  \sl For each nonzero fractional ideal $F$ of $D$,
    $F \subseteq F^{\ast_{\varphi}}$.  \rm

  \smallskip

  {\bf Claim 2}.  \sl For each nonzero $z \in L$,
  $(zD)^{\ast_{\varphi} }= zD$ (in particular, $D^{\ast_{\varphi}}
  =D$).  \rm

  Note that
  $$
  \begin{array}{rl}
      (zD)^{\ast_{\varphi} } &= \cap
      \left\{
          y\varphi\left(
                     \left(
                         \varphi^{-1}\left(
                                          y^{-1}zD
                                      \right)
                      \right)^{\ast}
                  \right)
      \mid \, y \in  L\,,\;  zD \subseteq yD \right\}  \\
  &\subseteq
  z\varphi\left(\left(\varphi^{-1}\left(D\right)\right)^{\ast}\right) =
  z\varphi\left(R^{\ast}\right) = z\varphi\left(R\right) = zD\,.
  \end{array}
  $$
  The conclusion follows from Claim 1.

  \smallskip

  {\bf Claim 3}.  \sl For each nonzero fractional ideal $F$ of $D$ and
for
  each nonzero $z\in L$, $(zF)^{\ast_{\varphi}}= zF^{\ast_{\varphi}}$.
  \rm

  Given $0 \neq z\in L$, for each nonzero $y \in L$, set $w:=yz
  \in L$. Then
  $$
  \begin{array}{rl}
  F^{\ast_{\varphi}} &=\cap

\left\{y\varphi\left(\left(\varphi^{-1}\left(y^{-1}F\right)\right)^{\ast}\right)

  \mid \, y \in L\,,\;  F \subseteq yD \right\}\\ &=\cap

\left\{\frac{w}{z}\varphi\left(\left(\varphi^{-1}\left(\frac{z}{w}F\right)\right)^{\ast}\right)

  \mid \, w \in L\,,\;  F \subseteq \frac{w}{z}D \right\} \\ &=\cap

\left\{{z}^{-1}w\varphi\left(\left(\varphi^{-1}\left({w}^{-1}zF\right)\right)^{\ast}\right)

  \mid \, w \in L\,,\; zF \subseteq wD \right\} \\ &=
  {z}^{-1}\left(\cap

\left\{w\varphi\left(\left(\varphi^{-1}\left({w}^{-1}zF\right)\right)^{\ast}\right)

  \mid \, w \in L\,,\; zF \subseteq wD \right\}\right) \\ &=
  z^{-1}(zF)^{\ast_{\varphi}}\,.
  \end{array}
  $$
  Hence, we conclude that $(zF)^{\ast_{\varphi}}=
  zF^{\ast_{\varphi}}$.

  \smallskip

  {\bf Claim 4}.  \sl For each pair of nonzero fractional ideals
  $F_{1} \subseteq F_{2}$ of $D$,  $(F_{1})^{\ast_{\varphi}}
  \subseteq (F_{2})^{\ast_{\varphi}}$.  \rm

  Note that if $y \in L$ and $F_{2}\subseteq yD$, then obviously
  $F_{1}\subseteq yD$, therefore
  $$
  \begin{array}{rl}
  (F_{2})^{\ast_{\varphi}} &=\cap

\left\{y\varphi\left(\left(\varphi^{-1}\left(y^{-1}F_{2}\right)\right)^{\ast}\right)

  \mid \, y \in L\,,\;  F_{2} \subseteq yD \right\} \\
  &\supseteq\cap

\left\{y\varphi\left(\left(\varphi^{-1}\left(y^{-1}F_{1}\right)\right)^{\ast}\right)

  \mid \, y \in L\,,\;  F_{1} \subseteq yD \right\} =
  (F_{1})^{\ast_{\varphi}}\,.
  \end{array}
  $$

  \smallskip

  {\bf Claim 5}.  \sl For each nonzero fractional ideal $F$ of $D$,
   $(F^{\ast_{\varphi}})^{\ast_{\varphi}} =F^{\ast_{\varphi}}$. \rm

  Note that from Claim 1, 2 and 4, if $y$ is a nonzero element of
  $L$,  we have

  \centerline{$ F\subseteq yD \; \Leftrightarrow \;
F^{\ast_{\varphi}}\subseteq
  (yD)^{\ast_{\varphi}} =
  yD\,,$}

    \noindent therefore
  $$
  \begin{array}{rl}
      (F^{\ast_{\varphi}})^{\ast_{\varphi}} &=\cap

\left\{y\varphi\left(\left(\varphi^{-1}\left(y^{-1}F^{\ast_{\varphi}}\right)\right)^{\ast}\right)

      \mid \, y \in L\,,\;  F^{\ast_{\varphi}} \subseteq yD \right\}
      \\
  &=\cap

\left\{y\varphi\left(\left(\varphi^{-1}\left(y^{-1}F^{\ast_{\varphi}}\right)\right)^{\ast}\right)

  \mid \, y \in L\,,\;  F \subseteq yD \right\}\,.
  \end{array}
  $$
  On the other hand,

  \centerline{$ F \subseteq yD \, \Rightarrow \, F^{\ast_{\varphi}}
\subseteq
y\varphi\left(\left(\varphi^{-1}\left(y^{-1}F\right)\right)^{\ast}\right)
  \ \Rightarrow \,
  y^{-1}F^{\ast_{\varphi}} \subseteq
\varphi\left(\left(\varphi^{-1}\left(y^{-1}F\right)\right)^{\ast}\right)$.}

   \noindent Therefore,

   \centerline{$ \varphi^{-1}\left(y^{-1}F^{\ast_{\varphi}}\right)
\subseteq
   \varphi^{-1} \left( \varphi \left( \left( \varphi^{-1} \left( y^{-1}F
   \right) \right)^{\ast} \right) \right) =
   \left(\varphi^{-1}\left(y^{-1}F\right)\right)^{\ast}\,,
   $}

   \noindent since $\left(\varphi^{-1}\left(y^{-1}F\right)\right)^{\ast}
   \supseteq \varphi^{-1}\left(y^{-1}F\right) \supseteq M = \mbox{\rm
   Ker}(\varphi)$.

   Now, we can conclude
   $$
   \begin{array}{rl}
      (F^{\ast_{\varphi}})^{\ast_{\varphi}} &=\cap

\left\{y\varphi\left(\left(\varphi^{-1}\left(y^{-1}F^{\ast_{\varphi}}\right)\right)^{\ast}\right)

      \mid \, y \in L\,,\;  F \subseteq yD \right\} \\
      &\subseteq
  \cap
  \left\{
         y\varphi \left( \left( \left( \varphi^{-1} \left( y^{-1}F
\right)
         \right)^{\ast} \right)^{\ast}\right) \mid \, y \in L\,,\; F
         \subseteq yD \right\} \\
  &=
  \cap
  \left\{
         y\varphi \left( \left( \varphi^{-1} \left( y^{-1}F \right)
         \right)^{\ast} \right) \mid \, y \in L\,,\; F \subseteq
         yD \right\}= F^{\ast_{\varphi}}\,,
  \end{array}
  $$
  and so, by Claim 1, $(F^{\ast_{\varphi}})^{\ast_{\varphi}}
  =F^{\ast_{\varphi}}$.  \hfill $\Box$

  \vspace{6pt}

  In case of a pullback of type \bf (\th$^{+}$) \rm the definition of
the
  star operation $\ast_{\varphi}$  given above is simplified as follows:

  \begin{prro}\label{prop:2.7} \it Let $T$, $K$, $M$, $k$, $D$,
$\varphi$,
  $L$, $S$ and $R$ be as in \rm \bf (\th$^{+}$)\it.  Let $\ast$ be a
given
  star operation on the integral domain $R$.  For each nonzero
fractional
  ideal $F$ of $D$, we have
      $$
      F^{\ast_{\varphi}} = \varphi\left(
\left(\varphi^{-1}(F)\right)^{\ast}
      \right) = \frac{\left(\varphi^{-1}(F)\right)^{\ast}}{M}\,.
      $$
      \end{prro}
       \rm

  \noindent {\bf Proof}.  For the extreme cases $M=(0)$ or $D=k$, it
trivially
  holds, so we may assume that $M\neq (0)$ and $D\subsetneq k$. We start
by proving the
  following:

  \bf Claim.  \sl Let $I$ be a fractional ideal of $R$ such that
  $M\subsetneq I\subseteq S={\varphi}^{-1}(L)$ and let $s\in S\setminus
  M$.  Then $(sI+M)^{\ast}=sI^{\ast}+M$.  \rm

  Choose $t\in S$ such that $st-1\in M$.  Then
  $t(sI+M)^{\ast}=(tsI+tM)^{\ast}\subseteq
  (tsI+M)^{\ast}=(I+M)^{\ast}=I^{\ast}$.  Therefore
  $st(sI+M)^{\ast}\subseteq sI^{\ast}$, so
  $st(sI+M)^{\ast}+M\subseteq sI^{\ast}+M\subseteq (sI+M)^{\ast}$.
  Put $m:=st-1$.  Since $m(sI+M)^{\ast}=(msI+mM)^{\ast}\subseteq
  M^{\ast}=M$ (where the last equality follows from the fact that $M$ is
a
  divisorial ideal of $R$), we have
  $st(sI+M)^{\ast}+M=(1+m)(sI+M)^{\ast}+M=(sI+M)^{\ast}$.  Thus we
  can conclude that $(sI+M)^{\ast}=sI^{\ast}+M$.

  \smallskip

  Now, let $F$ be a nonzero fractional ideal of $D$ and let
  $I:={\varphi}^{-1}(F)$.  For each element $y\in L$ such that
$F\subseteq
  yD$, we can find $s_y, t_y\in S\setminus M$ such that $\varphi
(s_y)=y$ and
  $\varphi (t_y)=y^{-1}$.  Using the above claim, we have
  $$
  \begin{array}{rl}
  F^{\ast_{\varphi}}
  &=\cap\{y\varphi(({\varphi}^{-1}(y^{-1}F))^\ast) \mid y\in L,\
  F\subseteq yD\} \\ &= \cap\{y\varphi ((t_yI+M)^\ast) \mid y\in L,\
  F \subseteq yD\}\\ &= \cap\{y\varphi (t_yI^\ast+M) \mid y\in L,\
  F\subseteq yD\} \\ &= \cap\{\varphi(s_y(t_yI^\ast+M)) \mid y\in
  L,\ F\subseteq yD\} \\ &= \cap\{\varphi(s_yt_yI^\ast+s_yM) \mid
  y\in L,\ F\subseteq yD\} \\ &= \cap\{\varphi(s_yt_yI^\ast+s_yM+M)
  \mid y\in L,\ F\subseteq yD\} \\ &= \cap\{\varphi(s_yt_yI^\ast+M)
  \mid y\in L,\ F\subseteq yD\} \\ &=
  \cap\{\varphi((s_yt_yI+M)^\ast)\mid y\in L,\ F\subseteq yD\} \\
  &=\cap\{\varphi(I^\ast)\mid y\in L,\ F\subseteq yD\}
  =\varphi(I^\ast) =\frac{I^\ast}{M}=
  \frac{\left(\varphi^{-1}(F)\right)^{\ast}}{M}\,.
  \end{array}
  $$
  \vskip -15 pt
  \hfill $\Box$

  \medskip

  \begin{reem}\label{rm:2.8}\rm As a consequence of Proposition
  \ref{prop:2.7} (and in the situation described in that statement) we
have the
  following:

  \sl If $I$ is a nonzero fractional ideal of $R$ such that
  $I\subseteq S$ and $sI\subseteq R$ for some $s\in S\setminus M$,
  then $I^\ast\subseteq S$ for any star operation $\ast$ on $R$. \rm As
  a matter of fact, $I^\ast \subseteq I^\ast S=I^\ast(M+sS)=I^\ast
   M+sI^\ast S\subseteq (IM)^\ast+(sI)^\ast S\subseteq M^\ast+S=M+S=S$.
   \end{reem}

  \begin{prro}\label{prop:2.9} Let $T$, $K$, $M$,  $k$, $D$, $\varphi$,
$L$, $S$ and
  $R$ be as in \rm \bf (\th$^{+}$)\it.  Let $\star$ be a given star
  operation on the integral domain $D$, let $\ast:= \star^{\varphi}$
  be the star operation on $R$ associated to $\star$ (which is defined
in
  Corollary \ref{cr:2.4}) and let $\ast_{\varphi} \; (=
  (\star^{\varphi})_{\varphi})$ be the star operation on $D$ associated
to
  $\ast$ (which is defined in Proposition \ref{pr:2.6}).  Then $\star
  =\ast_{\varphi} \ (= (\star^{\varphi})_{\varphi})$.
  \end{prro}

  \noindent {\bf Proof}.  For each nonzero fractional ideal $F$ of $D$
  and for each $y \in L$ such that $F \subseteq yD$, $J:= y^{-1}F$
  is a nonzero integral ideal of $D$.  Set $I_{y}:= \varphi^{-1}(J)
  =\varphi^{-1}(y^{-1}F) \,(\subseteq R)$.
  Note that $I_{y} $ is a nonzero ideal of $R$
  such that $M \subset I_{y}\subseteq R$, and so
  $\varphi(I_{y}) = I_{y}/M = J \, (\subseteq D) $.
  Moreover, we have
  $$
  \begin{array}{rl}
  (I_{y})^\ast &=\cap \left\{
                                x^{-1}\varphi^{-1} \left( \left(
\frac{xI_{y}+M}{M} \right)^{\star} \right)
                                \mid \, x\in I_{y}^{-1}\,\; x \neq
0\right\}\\
                                              &=\cap \left\{
                                x\varphi^{-1} \left( \left(
\frac{x^{-1}I_{y}+M}{M} \right)^{\star} \right)
                                \mid \; I_{y}\subseteq xR \subseteq K
\,\right\}  \\
       &= \left(\cap \left\{ x\varphi^{-1} \left( \left(
\frac{x^{-1}I_{y}+M}{M}
       \right)^{\star} \right) \mid \; I_{y}\subseteq
       xM\,, \; x \in K \right\}\right) \ \\
        & \mbox{ \;\;\,} \cap \left(\cap \left\{ x\varphi^{-1} \left(
\left(
        \frac{x^{-1}I_{y}+M}{M} \right)^{\star} \right) \mid \;
I_{y}\subseteq
        xR \subseteq K\,, \; \mbox{but } I_{y}\nsubseteq xM
\,\right\}\right)
        \\
         &= \left(\cap \left\{ xM \mid \; I_{y}\subseteq
       xM\,, \; x \in K \,\right\}\right) \ \\
        & \mbox{ \;\;\,}\cap  \left(\cap \left\{ x\varphi^{-1} \left(
\left(
        \frac{x^{-1}I_{y}+M}{M} \right)^{\star} \right) \mid \;
I_{y}\subseteq
        xR \subseteq K\,, \;\mbox{but } I_{y}\nsubseteq xM
\,\right\}\right)\,.
        \end{array}
  $$

  $\bullet$ \;\; For the first component of the previous
  intersection, note that since $M$ is maximal in $S$ and $M
  \subset I_{y} \subseteq R$, $I_{y}S =S$. On the other hand, $I_{y}
  \subseteq xM$, thus $ \varphi^{-1}(D) = R \subseteq S =I_{y}S
  \subseteq xMS =xM$. Therefore,  we have

  \centerline{$ \cap \left\{ xM \mid \; I_{y}\subseteq xM\subseteq K
  \,\right\}\supseteq \varphi^{-1}(D) \supseteq
  \varphi^{-1}\left((y^{-1}F)^{\star}\right)\,.
       $}

        $\bullet$ \;\; For the second component of the previous
        intersection, note that

        \centerline{$ x^{-1} I_{y} \subseteq R \, \mbox{ and } \, M
\subset I_{y}
        \subseteq R \; \Rightarrow \; x^{-1} I_{y}S \subseteq S \,
\mbox{ and }
        \, I_{y}S = S \; \Rightarrow \; x^{-1} \in S\,.
  $}

  \noindent On the other hand, if $ I_{y} \not\subseteq xM \;\;
  (I_{y} \subseteq xR)$ and $ x^{-1} \in S$, then $x^{{-1}}\in S
  \setminus M$,  and so  $ \varphi(x^{{-1}}) \in \varphi(S \setminus
  M) = L \setminus\{0\}\,.$ Note also that $ ({x^{-1}I_{y}+M})/{M} =
  \varphi(x^{-1})({I_{y}}/{M})\,. $

  Set

  \centerline{$
  I'_{y} := \varphi^{-1}\left((y^{-1}F)^{\star}\right) \;(\ \supseteq
  \varphi^{-1}\left(y^{-1}F\right) =: I_{y}\ )\ ,$}

  \noindent hence $ {I'_{y}}/{M} = \left( y^{-1}F\right)^{\star}=
  \left({I_{y}}/{M}\right)^{\star}\,.$

  Then we have
  $$
  \begin{array}{rl}
  &\cap \left\{ x\varphi^{-1} \left( \left(
        \frac{x^{-1}I_{y}+M}{M} \right)^{\star} \right) \mid \;
I_{y}\subseteq
        xR \subseteq K\,,\; \mbox{but } I_{y}\nsubseteq xM \,\right\}
        \\
  &\; = \cap \left\{ x\varphi^{-1} \left( \left(
        \varphi(x^{{-1}})\frac{I_{y}}{M} \right)^{\star} \right) \mid \;
I_{y}\subseteq
        xR \subseteq K\,,\; \mbox{but } I_{y}\nsubseteq xM \,\right\}
\\
  &\; = \cap \left\{ x\varphi^{-1} \left(\varphi(x^{{-1}}) \left(
        \frac{I_{y}}{M} \right)^{\star} \right) \mid \; I_{y}\subseteq
        xR \subseteq K\,,\; \mbox{but } I_{y}\nsubseteq xM \,\right\}
\\
       &\; = \cap \left\{ x\varphi^{-1} \left(\varphi(x^{{-1}})
\frac{I'_{y}}{M} \right) \mid \; I_{y}\subseteq
        xR \subseteq K\,,\; \mbox{but } I_{y}\nsubseteq xM \,\right\} \\

        &\; = \cap \left\{ x\left(x^{{-1}}I'_{y}+M \right) \mid \;
I_{y}\subseteq
        xR \subseteq K\,,\; \mbox{but } I_{y}\nsubseteq xM \,\right\}
\\
        &\; =  \cap \left\{ I'_{y}+xM  \mid \; I_{y}\subseteq
        xR \subseteq K\,,\; \mbox{but } I_{y}\nsubseteq xM \,\right\}  =

        I'_{y} = \varphi^{-1}\left((y^{-1}F)^{\star}\right)\,,
        \end{array}
        $$
        since for $x=1$ we have $I_y\subseteq xR\subseteq K$ but
        $I_y\not\subseteq xM$.

     \vspace{5 pt} Note that the first component of the intersection
representing $
     (I_{y})^\ast$ might not appear, but the second component
     necessarily appears, since at least for $x:=1$ we have that
     $I_{y}\subseteq xR \subseteq K\, \;\mbox{but } I_{y}\nsubseteq
xM$.  \
     Putting together the previous information about the two components
of
     the intersection, we have

        \centerline{$ \left(\varphi^{-1}(y^{-1}F)\right)^{\ast} =
(I_{y})^\ast =
        \varphi^{-1}\left((y^{-1}F)^{\star}\right)\,.
        $}

  \noindent Therefore we conclude that
        $$
      \begin{array}{rl}
      F^{\ast_{\varphi}} &= \cap

\left\{y\varphi\left(\left(\varphi^{-1}\left(y^{-1}F\right)\right)^{\ast}\right)

  \mid \, y \in L\,,\; F \subseteq yD \right\} \\ &= \cap \left\{
          y\varphi\left(\left(
                      I_{y}
                   \right)^{\ast}\right)
   \mid \, y \in L\,,\;  F \subseteq yD \right\} \\

   &= \cap
  \left\{
          y\varphi\left(
\varphi^{-1}\left((y^{-1}F)^{\star}\right)\right)
   \mid \, y \in L\,,\;  F \subseteq yD \right\} \\
    &= \cap
  \left\{
          y(y^{-1}F)^{\star}
   \mid \, y \in L\,,\; F \subseteq yD \right\} \\
     &= \cap
  \left\{
          yy^{-1}F^{\star}
   \mid \, y \in L\,,\; F \subseteq yD \right\}
   =F^{\star}\,.  \;\;\;\; \; \hskip 3.7cm \Box
  \end{array}
  $$

  \begin{reem}\label{rmk:2.10} \sl With the notation and hypotheses of
  Proposition \ref{prop:2.9}, for each nonzero fractional ideal $F$ of
  $D$, we have

  \centerline{$ F^\star
=\varphi\left(\varphi^{-1}(F)^{\star^{\varphi}}\right)\,.
  $}

  \noindent\rm  As a matter of fact, by the previous proof and
Proposition
  \ref{prop:2.7}, we have that $F^\star = F^{\ast_{\varphi}} =
  {\varphi^{-1}(F)^{\star^{\varphi}}}/{M}$.
  \end{reem}

  \vskip 6pt

  \begin{coor}\label{coro:2.10} Let $T$, $K$, $M$,  $k$, $D$, $\varphi$,
$L$, $S$ and
  $R$ be as in \rm \bf (\th$^{+}$)\it.

  \balf
  \bf \item \it The map  \ \bf $(-)_{\varphi} :$ Star\rm$(R)
\rightarrow$ \bf
  Star\rm$(D)$,\ $\ast \mapsto \ast_{\varphi}$, \it \ is
order--preserving
  and surjective.

  \bf \item \it The map \ \bf $(-)^{{\varphi}} :$ \bf Star\rm$(D)
  \rightarrow$ \bf Star\rm$(R)$,\ $\star \mapsto \star^{\varphi}$, \it \
  is order--preserving and injective.

  \bf \item \it Let $\star$ be a star operation on $D$. Then  for each
nonzero
  ideal $I$ of $R$ with $M \subset I \subseteq R$,

  \centerline{$ I^{\star^{\varphi}} =
  \varphi^{-1}\left(\left(\varphi(I)\right)^\star\right)\,.$}
  \ealf
  \end{coor}

  \noindent {\bf Proof}.  \bf (a) \rm and \bf (b) \rm are
straightforward
  consequences of the definitions and  Proposition \ref{prop:2.9}, since
  $(-)^{{\varphi}} $ is a right inverse of $(-){_{\varphi}}$ (i.e.,
  $(-)_{\varphi}\circ(-)^{{\varphi}}=$ {\bf 1}$_{\mbox{\tiny \bf Star\rm
  }(D)}$).

  \bf (c) \rm Let $\ast :=\star^{\varphi}$.  Then by Proposition
  \ref{prop:2.9}, we know that $\ast_{\varphi }= \star$.  Therefore,
using
  Proposition \ref {prop:2.7}, we have

  \centerline{$ \left(\varphi(I)\right)^\star =
  \left(\varphi(I)\right)^{\ast_{\varphi }} =
  \frac{\left(\varphi^{-1}(\varphi(I))\right)^{\ast}}{M}=
  \frac{I^\ast}{M}= \frac{I^{\star^{\varphi}}}{M}\,,
  $}

  \noindent and hence $\varphi^{-1}\left(\left(\varphi(I)\right)^\star
\right) =
  I^{\star^{\varphi}}$.

  \vskip -15pt \hfill $\Box$

  \vspace{12pt}

  The next result shows how the composition map

  \centerline{$(-)^{{\varphi}} \circ\ (-)_{\varphi} :$\ \bf Star\rm$(R)
  \rightarrow $ \bf Star\rm$(R) $}

  \noindent compares with the identity map.

  \begin{thee}\label{thr:2.11} Let $T$, $K$, $M$, $k$, $D$, $\varphi$,
  $L$, $S$ and $R$ be as in \rm \bf (\th$^{+}$)\it.  Assume that
  $D\subsetneq k$.  Then for each star operation $\ast$ on $R$,

  \centerline{$ \ast \leq \left((\ast)_{\varphi}\right)^{\varphi}\,.
  $}
  \end{thee}

  \noindent \bf Proof. \rm
  Let $I$ be a nonzero integral ideal of $R$.  For each nonzero $x\in
I^{-1}$, if $xI \not\subseteq M$, then by Proposition \ref{prop:2.7},
$\left(\frac{xI+M}{M}\right)^{\ast_{\varphi}} =
\frac{(xI+M)^{\ast}}{M}\supseteq \frac{(xI)^{\ast}+M}{M}$.
Now using the fact $M^{\ast}=M$ for $M\neq (0)$, we have
$$
  \begin{array}{rl}
      I^{(\ast_{\varphi})^{\varphi}} &=
      \cap \left\{

x^{-1}\varphi^{-1}\left(\left(\frac{xI+M}{M}\right)^{\ast_{\varphi}}\right)

  \mid \, x \in I^{-1}\,, \; x \neq 0 \right\} \\
  &= \left(\cap \left\{

x^{-1}\varphi^{-1}\left(\left(\frac{xI+M}{M}\right)^{\ast_{\varphi}}\right)

  \mid \, x \in I^{-1}\,, \; x \neq 0\,, \; xI \not\subseteq M
\right\}\right)
  \cap \\
  & \mbox{ } \hskip 0.8cm \left(\cap \left\{ x^{-1}M \mid
  \, x \in I^{-1}\,, \; x \neq 0\,, \; xI \subseteq M \right\}\right) \\

  &\supseteq \left(\cap \left\{
  x^{-1}\varphi^{-1}\left(\frac{(xI)^\ast+M}{M}\right) \mid \, x \in
  I^{-1}\,, \; x \neq 0\,, \; xI \not\subseteq M \right\} \right) \cap
\\
 & \mbox{ } \hskip 0.8cm \left(\cap \left\{ x^{-1}M^\ast \mid \, x \in
I^{-1}\,,
 \; x \neq 0\,, \; I \subseteq x^{-1}M \right\}\right) \\
 &\supseteq \left(\cap \left\{ x^{-1}\left({(xI)^\ast+M}\right) \mid \,
x
 \in I^{-1}\,, \; x \neq 0\,, \; xI \not\subseteq M \right\}\right) \cap

 I^\ast \\
 &\supseteq \left(\cap \left\{ x^{-1}\left({(xI)^\ast}\right) \mid \, x
\in
 I^{-1}\,, \; x \neq 0\,, \; xI \not\subseteq M \right\}\right) \cap
 I^\ast  = I^\ast \,.
      \end{array}
      $$
      \vskip -15 pt
  \hfill $\Box$\
      \vskip 5pt

      In Section 3, we will show that in general  $
  \ast \lneq \left((\ast)_{\varphi}\right)^{\varphi}$. However,
  in some relevant cases,  the inequality
    is, in fact, an equality:

  \begin{coor}\label{coro:2.13} Let $T$, $K$, $M$,  $k$, $D$, $\varphi$,
$L$, $S$ and
  $R$ be as in Theorem
  \ref{thr:2.11}.  Then

  \centerline{$ v_{R} =
  ((v_{R})_{\varphi})^{\varphi}\,;\;\;\;\;\;\;\;
  (v_{D})^{\varphi}
  =v_{R}\,;\;\;\;\;\;\;\;  (v_{R})_{\varphi} =v_{D}\,.
  $}
  \end{coor}

  \noindent {\bf Proof}.  Use
  Proposition \ref{prop:2.9}, Corollary \ref{coro:2.10} (b), Theorem
  \ref{thr:2.11} and \cite[Theorem 34.1 (4)]{G}.  More precisely, note
that
  $(v_{R})_{\varphi} \leq v_{D}$, and so $v_{R} \leq
  ((v_{R})_{\varphi})^{\varphi} \leq (v_{D})^{\varphi}\leq v_{R}$.  On
the
  other hand, if $(v_{R})_{\varphi} \lneq v_{D}$, then
  $v_{R}= \left((v_{R})_{\varphi}\right)^{\varphi} \lneq
  (v_{D})^{\varphi}$, which is a contradiction.  \hfill $\Box$

  \vskip 10pt
  Our next goal is to apply the previous results for giving a
componentwise
  description of the ``pullback'' star operation $\diamond$ considered
in
  Proposition \ref{pr:2.1}.

  \vskip 5pt

  \begin{prro}\label{thr:2.12}
  Let $T$, $K$, $M$,  $k$, $D$, $\varphi$, $L$, $S$ and
  $R$ be as in \rm \bf (\th$^{+}$)\it.  Assume that $M\neq (0)$ and
  $D\subsetneq k$.  Let

  \centerline{$\Phi : \mbox{\bf Star}(D)\times \mbox{\bf Star}(T)
\rightarrow
  \mbox{\bf Star}(R)$\,, \ $({\star}_D, {\star}_T)\mapsto \diamond :=
  ({\star}_D)^{\varphi}\wedge ({\star}_T)^{\iota}$,}

  \noindent be the map considered in Proposition \ref{pr:2.1} and
Corollary
  \ref{cr:2.5}.  The following properties hold:

  \begin{enumerate}
  \bf \item[(a)] \it ${\diamond}_{\varphi}={\star}_D$.

  \bf \item[(b)] \it ${\diamond}_{\iota}=(v_R)_{\iota}\wedge {\star}_T \
  (\in \mbox{\bf Star}(T))$.

  \bf \item[(c)] \it $\diamond=({\diamond}_{\varphi})^{\varphi}\wedge
  ({\diamond}_{\iota})^{\iota}$.
  \end{enumerate}

  \end{prro}

  \noindent {\bf Proof}.  \bf (a) \rm Without loss of generality, we
  only consider the case of integral ideals of $D$.  Let $J$ be a
nonzero
  integral ideal of $D$ and let $I:={\varphi}^{-1}(J)$.  Since
  $M\subsetneq I\subseteq R$, we have $IS=S$, where
  $S:={\varphi}^{-1}(L)$, and so $IT=T$.  Therefore, by Proposition
  \ref{prop:2.7} and Corollary \ref{coro:2.10} (c),
$J^{{\diamond}_{\varphi}}=\varphi(I^{\diamond})=\varphi(I^{({\star}_D)^{\varphi}}\cap
  I^{({\star}_T)^{\iota}})=\varphi(I^{({\star}_D)^{\varphi}}\cap
  (IT)^{{\star}_T})=\varphi(I^{({\star}_D)^{\varphi}}\cap
T)=\varphi(I^{({\star}_D)^{\varphi}})=\varphi({\varphi}^{-1}(J^{{\star}_D}))=J^{{\star}_D}$.

  \bf (b) \rm Without loss of generality, we only consider  the case of
  integral ideals of $T$.  Let $I$ be a nonzero ideal of $T$
   (in particular, $I$ is a fractional ideal of $R$).  Then
for
   each $x\in I^{-1}=(R:_{K} I)$, we have $xIT=xI\subseteq R$, so
   $xI\subseteq (R:_{K}T)=M$.  Therefore
  $$
  \begin{array}{rl}
      I^{({\star}_D)^{\varphi}}
  &=\bigcap\{x^{-1}{\varphi}^{-1}((\varphi(xI))^{{\star}_D}) \,|\,
  x\in I^{-1},\ x \neq 0 \}\\ &=  \bigcap\{x^{-1}M \,|\, x\in
  I^{-1},\ x \neq 0 \}=I^{v_R}\,,
  \end{array}
  $$
  and so

  \centerline{$
  I^{{\diamond}_{\iota}}=I^{\diamond}=I^{({\star}_D)^{\varphi}}\cap
  I^{{\star}_T}=I^{v_R}\cap I^{{\star}_T}= I^{(v_R)_{\iota}}\cap
  I^{{\star}_T}=I^{(v_R)_{\iota}\wedge {\star}_T}\,.
  $}

  \noindent Note that $I^{{\diamond}_{\iota}} \ (\subseteq I^{v_R})$
  is an ideal of $R$.  Moreover, $I^{{\diamond}_{\iota}}$ is an
  ideal of $T$, because for each nonzero $x\in T$,
  $xI^{{\diamond}_{\iota}}=x(I^{v_R}\cap
  I^{{\star}_T})=(xI)^{v_R}\cap (xI)^{{\star}_T}\subseteq
  I^{v_R}\cap I^{{\star}_T}=I^{{\diamond}_{\iota}}$.  Finally, since
  $\star_{T}$ is a star operation on $T$, it is easy to check that
  ${\diamond}_{\iota}$ (restricted to $\boldsymbol{F}(T)$) belongs
  to $\mbox{\bf Star}(T)$.

  \bf (c) \rm Since $\diamond\leq v_R\leq ((v_R)_{\iota})^{\iota}$,
   (using also Example \ref{ex:1.5}) we have that $$
  \begin{array}{rl} \diamond \hskip -0.2 pt= \hskip -0.2
  pt ({\star}_D)^{\varphi}\wedge ({\star}_T)^{\iota} \hskip -0.2 pt
  = & \hskip -6 pt  ({\star}_D)^{\varphi}\wedge
  ((v_R)_{\iota})^{\iota}\wedge ({\star}_T)^{\iota}\hskip -0.2 pt \\
  =& \hskip -6 pt ({\star}_D)^{\varphi}\wedge ((v_R)_{\iota}\wedge
  {\star}_T)^{\iota} \hskip -0.2 pt = \hskip -0.2 pt
  ({\diamond}_{\varphi})^{\varphi}\wedge
  ({\diamond}_{\iota})^{\iota}.
  \end{array}
  $$

  \vskip -20pt \hfill $\Box$

  \begin{exxe} \label{ex:2.15} \sl With the same notation and hypotheses
  of Proposition \ref{thr:2.12}, we show that, in general,
  ${\diamond}_{\iota}\neq {\star}_T$.  \rm

  \bf (1) \rm Let $T:=k[X,Y]_{(X,Y)}$ and let $M:=(X,Y)T$.  Then $T$ is
a
  2-dimensional  local UFD. Choose a subfield $D:=L$ of $k$ such that
  $[k:L]=2$.  In this situation we have that $T\subseteq
(R:_{K}M)\subseteq
  (T:_{K}M)$, and $(T:_{K}M)=T$ because $T$ is 2-dimensional local UFD
(hence, Krull)
  with maximal ideal $M$.  Therefore,  $T=(R:_{K}M)$.  By Example
  \ref{ex:2.3} (d), if ${\star}_T:=v_T$, then $\diamond=v_R$ and
$M^{v_T} =
  T$.  But $M^{{\diamond}_{\iota}}=M^{\diamond}=M^{v_R}=M\neq
  T=M^{v_T}=M^{{\star}_T}$.

  \bf (2) \rm Note that ${\diamond}_{\iota}\neq {\star}_T$, even if
$L=k$.
  It is sufficient to consider a slight modification of the previous
example.
  Let $D$ be any integral domain (not a field) with quotient field $L$.
Let
  $T:={L[X, Y]}_{(X, Y)}$ and let $M:=(X, Y)T$.
  Set
  $\overline{\diamond}:=(v_{D})^{\varphi}\wedge (v_{T})^{\iota}$.  Then
  ${M}^{\overline{{\diamond}}_{\iota}}={M}^{\overline{\diamond}}=
  {M}^{(v_{D})^{\varphi}} \cap {M}^{(v_{T})^{\iota}}= {M}^{v_{R}} \cap
  {M}^{(v_{T})^{\iota}} = {M}$, because ${M}^{v_{R}}= {M}$ and
  ${M}^{(v_{T})^{\iota}}= (MT)^{v_{T}} ={M}^{v_{T}}=T$.

  \end{exxe}

  \begin{reem} \rm \bf (a) \rm Note that, with the same notation and
  hypotheses of Proposition \ref{thr:2.12}, \sl the map $\Phi$ is not
  one-to-one in general.  \rm

  This fact follows immediately  from Example \ref{ex:2.15} and
Proposition
  \ref{thr:2.12} (b) and (c), since

  \centerline{$({\star}_D)^{\varphi}\wedge ({\star}_T)^{\iota} =
\diamond=({\diamond}_{\varphi})^{\varphi}\wedge
  ({\diamond}_{\iota})^{\iota}$\,.}

  \bf (b) \rm In the same setting as above, \sl the map $\Phi$ is not
onto in
  general. \rm

  For instance, in the situation described in Example \ref{ex:2.3}
  (d), we have that $d_R\not\in \mbox{Im}(\Phi)$.
   Another example, even in case $L=k$, is given next.

  \end{reem}

  \begin{exxe} \label{ex:2.17} \sl Let $D$ be a 1-dimensional discrete
  valuation domain with quotient field $L$.  Set $T:=L[X^2, X^3]$,
  $M:=X^2L[X]=XL[X] \cap T$ and $K:=L(X)$.  Let $\varphi$ and $R$
  be as in \rm \bf (\th$^{+}$)\rm .  \sl Then $v_{R}\notin$ \rm
Im$(\Phi)$.\rm

  Note that, for each $\diamond\in$ Im$(\Phi)$,
  $\diamond\leq (v_D)^{\varphi}\wedge (v_T)^{\iota}\leq v_R$.  In order
to
  show that $v_R\not\in$ Im$(\Phi)$, it suffices to prove that
  $(v_D)^{\varphi}\wedge (v_T)^{\iota}\neq v_R$.  The fractional
overring $T$ of $R$
  is not a divisorial ideal of $R$, since
  $T^{v_{R}}=(R:_{K}(R:_{K}T))=(R:_{K}M)\supseteq L[X]\supsetneq T$.
Therefore,
  $T^{(v_D)^{\varphi}\wedge (v_T)^{\iota}}=T^{v_R\wedge
  (v_T)^{\iota}}=T^{v_R}\cap T^{(v_T)^{\iota}}=T^{v_R}\cap
  T^{v_T}=T^{v_R}\cap T=T\subsetneq T^{v_R}$.

  \end{exxe}

  \begin{thee}\label{thr:2.18} With the notation and hypotheses of
Proposition
  \ref{thr:2.12}, set

   \centerline{$\mbox{\bf Star}(T; v_{R}):=\{{\star}_T\in \mbox{\bf
   Star}(T) \mid \,  {\star}_T\leq (v_R)_{\iota}\}$.}

   \noindent Then
  \begin{enumerate}
  \bf \item[(a)] \it $ \mbox{\bf Star}(T; v_{R}) = \{{\star}_T\in
  \mbox{\bf Star}(T) \mid \, \left(v_{R}\wedge
  ({\star}_T)^\iota\right)_{\iota} ={\star}_T \ \} \\
   {\, }\hskip 1.75cm =
   \{\ast_{\iota}
   \mid \, \ast\in \mbox{\bf Star}(R) \}\  \cap\ \mbox{\bf Star}(T) \\
  {\, }\hskip 1.75cm = \{\ast_{\iota} \mid \, \ast\in \mbox{\bf Star}(R)
  \,\mbox{ and }\; T^{\ast}=T\}\,.$

  \bf \item[(b)] \it The restriction $\Phi':= \Phi|_{\mbox{\footnotesize
{\bf
  Star}}(D) \times \mbox{\footnotesize {\bf Star}}(T; v_{R})}$ is
  one-to-one.

  \bf \item[(c)] \it $\mbox{\rm Im}(\Phi') = \mbox{\bf Star}(R;
\mbox{\bf
  (\th$^{+}$)}) := \{ \ast \in \mbox{\bf Star}(R) \mid T^\ast = T \mbox{

  and }\ \ast =({\ast}_{\varphi})^{\varphi}\wedge
({\ast}_{\iota})^{\iota}
  \}$.
  \end{enumerate}
  \end{thee}

  \noindent {\bf Proof}.  \bf (a) \rm We start by proving the following:

  \bf Claim.  \sl Let ${\star}_T \in \mbox{\bf Star}(T; v_{R})$ and
  let ${\star}_D \in \mbox{\bf Star}(D)$ be any star operation on $D$.
  Set, as usual,
  $\diamond :=({\star}_D)^{\varphi}\wedge ({\star}_T)^{\iota}$. \
  Then $\diamond_{\iota}= {\star}_T$.  \rm

  Note that, by Corollary \ref{coro:2.13},  $\diamond =
  \Phi((\star_D, {\star}_T)) \leq \overline{\diamond}:= \Phi((v_D,
  {\star}_T))= (v_D)^\varphi \wedge ({\star}_T)^{\iota} =v_R \wedge
  ({\star}_T)^{\iota}\in \mbox{\bf Star}(R)\,.$  Hence, by using
  Theorem \ref{thr:2.12} (b), Examples \ref{ex:1.2} (h) and
  \ref{ex:1.5} (b), we have $ (v_R)_{\iota}\wedge {\star}_T
  ={\diamond}_{\iota} \leq {\overline{\diamond}}_{\iota}=(v_R\wedge
  ({\star}_T)^{\iota})_{\iota}= (v_R)_{\iota}\wedge
  (({\star}_T)^{\iota})_{\iota} =(v_R)_{\iota} \wedge {\star}_T\,,$
  \ thus ${\diamond}_{\iota} = {\overline{\diamond}}_{\iota}=
  {\star}_T\,,$ \ because ${\star}_T \in \mbox{\bf Star}(T; v_{R})$.

  \smallskip

  From the previous argument we also deduce  that

  \centerline{ ${\star}_T\leq (v_R)_{\iota} \;\; \Leftrightarrow\;\;
  \left(v_{R}\wedge ({\star}_T)^\iota\right)_{\iota}
  ={\star}_T$\,.}

   Now, let $\ast\in \mbox{\bf Star}(R)$ be a star operation on
  $R$ such that $\ast_{\iota}\in \mbox{\bf Star}(T)$.  Then obviously
  $\ast_{\iota}\leq (v_R)_{\iota}$, whence $\ast_{\iota} \in \mbox{\bf
  Star}(T; v_{R})$, and $T^{\ast}=T^{\ast_{\iota}}=T$.

  If $\ast \in \mbox{\bf Star}(R)$ is such that $T^{\ast}=T$\,, \ then
  clearly we have $\ast_{\iota}\in \mbox{\bf Star}(T)$.

   If ${\star}_T \in \mbox{\bf Star}(T; v_{R})$, then by the Claim,
   ${\star}_T={\overline{\diamond}}_{\iota}$ with $\overline{\diamond}
\in
   \mbox{\bf Star}(R)$, hence ${\star}_T \in \{\ast_{\iota} \mid \,
\ast\in \mbox{\bf
   Star}(R) \} \cap \mbox{\bf Star}(T)$.

    \bf (b) \rm is a straightforward consequence of the Claim and of
    Proposition \ref{thr:2.12} (a).

   \bf (c) \rm follows from the Claim and from Proposition
\ref{thr:2.12}
   (a) and (c).  \hfill $\Box$

  \vskip 10 pt We next apply some of the theory developed above for
  answering  a problem posed by D. F. Anderson in 1992.

  \begin{exxe}\label{exe:2.14} \bf
(``$\boldsymbol{D+M}$''--constructions).\rm

    Let $T$ be an integral domain of the type $k+M$, where $M$ is a
maximal
    ideal of $T$ and $k$ is a subring of $T$ canonically isomorphic to
the
    field $T/M$, and let $D$ be a subring of $k$ with field of quotients
$L$\
    ($\subseteq k$).  Set $R:= D+M$.  Note that $R$ is a faithfully flat
    $D$--module.

    Given a star operation $\ast$ on $R$, D.F. Anderson \cite[page 835]
{A:1988}
    defined a star operation on $D$ in the following way: for each
nonzero
    fractional ideal $F$ of $D$, set

   \vskip -4pt \centerline{$ F^{\ast_{D}}:= (FR)^\ast \cap L \,.  $}

    \vskip 2pt \noindent Note that $FR = F+M$.  From \cite[Proposition
5.4 (b)]
    {A:1988} it is known that for each nonzero fractional ideal $F$ of
$D$,
    \bara \bf \item \sl $F^{\ast_{D}}+M = (F+M)^\ast$\,; \bf \item \sl
    $F^{\ast_{D}}= (F+M)^\ast \cap L = (F+M)^\ast \cap k$\,.  \eara

    \bf Claim.  \sl If $\varphi: R \rightarrow D$ is the canonical
    projection and if $\ast_{\varphi}$ is the star operation defined in
    Proposition \ref{pr:2.6}, then $\ast_{D}= \ast_{\varphi}$.

    In particular, by \rm \cite[Proposition 2 (a), (c)]{A:1992}\sl
    , we deduce that

    \balf
    \bf \item \sl $(d_{R})_{\varphi} = d_{D}\,,\;\;\;\;
(t_{R})_{\varphi} = t_{D}\,,\;\;\;\;
    (v_{R})_{\varphi} = v_{D}\,,\;$ and
   \bf \item \sl  $({\ast_f})_{\varphi} = (\ast_{\varphi})_{f}$\,.
   \ealf
   \rm

    Note that if $y$ is a nonzero element of the quotient field $L$ of
$D$, then
    $y$ belongs to $k$, and thus, $y$ is a unit in $T$ and so
$y^{-1}M=M$. Therefore, for
    each $ y \in L$ such that $ F \subseteq yD$, we have
    $$
    \begin{array}{rl}

y\varphi\left(\left(\varphi^{-1}\left(y^{-1}F\right)\right)^{\ast}\right)
&=
     y\varphi\left(\left(y^{-1}F+M\right)^{\ast}\right) =
     y\varphi\left(\left(y^{-1}F+y^{-1}M\right)^{\ast}\right)\\
     &= y\varphi\left(y^{-1}\left(F+M\right)^{\ast}\right)=
     y\varphi\left(y^{-1}(F^{\ast_{D}}+M) \right)\\
     &= y\varphi\left(y^{-1}F^{\ast_{D}}+y^{-1}M \right)=
     y\varphi\left(y^{-1}F^{\ast_{D}}+M \right)\\
     &= y \left(y^{-1}F^{\ast_{D}}\right)=F^{\ast_{D}}\,,
     \end{array}
     $$
     hence (Proposition \ref{pr:2.6}) $F^{\ast_{\varphi}}=
F^{\ast_{D}}$.

\smallskip

  By applying Proposition \ref{prop:2.9} and Corollary \ref{coro:2.10}
  (a) to the particular case of $R =D+M $ (special case of  \bf
(\th$^{+}$)\rm ),
  we know that the map

  \centerline{\, \bf $(-)_{\varphi} :$ Star\rm$(D+M) \rightarrow$ \bf
  Star\rm$(D)$,\, \; $\ast \mapsto \ast_{\varphi}=\ast_{D}$ \,,}

  \noindent is surjective and order-preserving and it has the injective
  order-preserving map

  \centerline{\, $(-)^{_{\varphi}} :$ \bf Star\rm$(D) \rightarrow$ \bf
  Star\rm$(D+M)$, \,\; $\star \mapsto \star^{\varphi}$ \,,}

  \noindent as a right inverse.  This fact gives a complete positive
answer to a
  problem posed by D.F. Anderson (cf.  \cite [page 226]{A:1992}).

    \end{exxe}


  \section{Transfer of star properties}

      In this section we want to investigate the general problem of the
      transfer --in the
      pullback setting-- of some relevant properties concerning the star
      operations involved.  In particular, we pursue the work initiated
by D.F.
      Anderson in \cite{A:1992} for the case of the
``$D+M$''--constructions. We
      start by studying which of the properties (a) and (b) of Example
      \ref{exe:2.14} hold in a more
      general setting.

    \vskip 10pt

  \begin{prro}\label{prop:3.1} Let $T$, $K$, $M$,  $k$, $D$, L,
$\varphi$ and $R$ be
  as in \rm \bf (\th$^{+}$)\it.  \balf

  \bf \item \it Let ${\cal{R}}:= \{ R_{\lambda} \mid\, \lambda \in
  \Lambda\}$ be a family of overrings of $R$ contained in $T$ such that
  $\cap \{ R_{\lambda}\mid\, \lambda \in \Lambda\} = R$, and let
  ${\cal{D}}:= \{ D_{\lambda}:= \varphi(R_{\lambda}) \mid\, \lambda \in
  \Lambda\}$ be the corresponding family of subrings of $k$ (with $\cap
\{
  D_{\lambda} \mid\, \lambda \in \Lambda\} = D$), then

  \centerline{$ (\star_{\cal{R}})_{\varphi} = \star_{\cal{D}}\,.$}

  \vskip -6pt \bf \item \it If ${\cal{D}}:= \{ D_{\lambda} \mid\,
\lambda
  \in \Lambda\}$ is a family of overrings of $D$ such that $\cap \{
  D_{\lambda}\mid\, \lambda \in \Lambda\} = D$ and if ${\cal{R}}:= \{
  R_{\lambda}:= \varphi^{-1}(D_{\lambda}) \mid\, \lambda \in \Lambda\}$
is
  the corresponding family of subrings of $T$ (with $\cap \{ R_{\lambda}
  \mid\, \lambda \in \Lambda\} = R$), then in general

  \centerline{$
  \star_{\cal{R}} \leq (\star_{\cal{D}})^{\varphi}\,.$}
  \ealf
  \end{prro}

  \noindent {\bf Proof}.  \rm \bf (a) \rm Note that in the present
situation
  $\varphi^{-1}(D_{\lambda}) = R_{\lambda}$ for each $\lambda \in
  \Lambda$\,, \ $ D = \cap \{ D_{\lambda} \mid\, \lambda \in \Lambda\}$,\
and
  for each nonzero fractional ideal $J$ of $D$,
  $J^{(\star_{\cal{R}})_{\varphi}} = \varphi\left(
  (\varphi^{-1}(J))^{\star_{\cal{R}}}\right)$\ (Proposition
  \ref{prop:2.7}).  Moreover,
  $$\begin{array}{rl}
   \varphi\left(
  (\varphi^{-1}(J))^{\star_{\cal{R}}}\right)&= \varphi\left(\cap\{
  \varphi^{-1}(J)R_{\lambda} \mid \, \lambda \in \Lambda\}\right) \\
  &= \varphi\left(\cap\{ \varphi^{-1}(J)\varphi^{-1}(D_{\lambda})
  \mid \, \lambda \in \Lambda\}\right)\\ &=
  \varphi\left(\varphi^{-1}(\cap\{ JD_{\lambda}\mid \, \lambda \in
  \Lambda\})\right)=
  \varphi\left(\varphi^{-1}(J^{\star_{\cal{D}}})\right)=
  J^{\star_{\cal{D}}}\,.
  \end{array}
  $$

  \bf (b) \rm Note that
  $\varphi(R_{\lambda})=\varphi(\varphi^{-1}(D_{\lambda})) =
D_{\lambda}$
  for each $\lambda \in \Lambda$.  Therefore, by (a), $
  ({\star_{\cal{R}}})_{\varphi}= \star_{\cal{D}}$, thus $
  (({\star_{\cal{R}}})_{\varphi})^{\varphi}=
  ({\star_{\cal{D}}})^{\varphi}$.  If $D=k$, then $D=L$ is a field, thus
  ${\cal D}=\{D\}$ and ${\cal R}=\{R\}$.  So, obviously,
  $\star_{\cal{R}}=d_R\leq ({\star_{\cal{D}}})^{\varphi}$.  If
  $D\subsetneq k$, then the conclusion follows from Theorem
  \ref{thr:2.11}.\hfill $\Box$

  \vskip 10pt Proposition \ref{prop:3.1} (a) can be
  generalized to a statement concerning a surjective homomorphism
between two
  integral domains:

  \begin{prro} \label{propo:3.2} \it Let $R$, $K$, $M$, $D$, $\varphi$
be
  as in Corollary \ref{cr:2.4}.  Let $\{\ast_{\lambda} \mid \lambda \in
  \Lambda \}$ be a family of star operations of $R$.  Then

  \centerline{$(\wedge {\ast}_{\lambda})_{\varphi}=\wedge
  ({\ast}_{\lambda})_{\varphi}$\,.}
  \end{prro}

  \noindent {\bf Proof}.  Let $J$ be a nonzero fractional ideal of $D$
and
  let $y$ be in the quotient field $L$ of $D$. Then
  $$
  \begin{array}{rl}
  J^{\wedge ({\ast}_{\lambda})_{\varphi}} & ={\cap} \{
  J^{({\ast}_{\lambda})_{\varphi}}\mid \lambda \in \Lambda \}
  \\ & = {\cap} \left\{
   (\cap\{y\varphi(({\varphi}^{-1}(y^{-1}J))^{{\ast}_{\lambda}})\mid
J\subseteq
   yD\}) \, \mid \lambda \in \Lambda \right\} \\
  & =\cap \left\{ y({\cap} \{
  \varphi\left(({\varphi}^{-1}(y^{-1}J))^{{\ast}_{\lambda}}\right)
  \mid \lambda \in \Lambda \}) \, \mid J\subseteq yD \right\} \\ &=
  \cap \left\{y\varphi\left({\cap} \{
  ({\varphi}^{-1}(y^{-1}J))^{{\ast}_{\lambda}}\mid \lambda \in
  \Lambda \}\right) \, \mid J\subseteq yD \right\} \\
   &
  =\cap \left\{y\varphi\left( ({\varphi}^{-1}(y^{-1}J))^{\wedge
  {\ast}_{\lambda}}\right) \, \mid J\subseteq yD
  \right\} = J^{(\wedge {\ast}_{\lambda})_{\varphi}}\,.
  \end{array}
  $$
  \vskip -16 pt\hfill $\Box$

  \vskip 15 pt

  \begin{prro} \label{prop:3.2} \it Let $R$, $K$, $M$, $D$, $\varphi$ be
as in
  Corollary \ref{cr:2.4}. Then

  \centerline{$
  (d_{R})_{\varphi}= d_{D}\,.$}
  \end{prro}

  \noindent \bf Proof.  \rm For each nonzero fractional ideal $J$ of
$D$, we
  have $$
  \begin{array}{rl} J^{(d_{R})_{\varphi}}&= \cap

\left\{y^{-1}\varphi\left(\left(\varphi^{-1}\left(yJ\right)\right)^{d_{R}}\right)

  \mid \, y \in J^{-1}\,,\; y \neq 0 \right\} \\ &= \cap

\left\{y^{-1}\varphi\left(\left(\varphi^{-1}\left(yJ\right)\right)\right)

  \mid \, y \in J^{-1}\,,\; y \neq 0 \right\}\\ &= \cap
  \left\{y^{-1}\left(yJ\right) \mid \, y \in J^{-1}\,,\; y \neq 0
  \right\}= J = J^{d_{D}}\,.\\
  \end{array}
  $$
  \vskip -15pt \hfill $\Box$

  \vskip 10pt

      The next couple of examples explicitly show that the inequalities
in
      Theorem \ref{thr:2.11} and Proposition \ref{prop:3.1} (b) can be
strict
      inequalities (i.e., $ \ast \lneq
\left((\ast)_{\varphi}\right)^{\varphi}$
      and $\star_{\cal{R}} \lneq (\star_{\cal{D}})^{\varphi}$).

  \begin{exxe} \label{ex:3.3}
   \sl Let $T$, $K$, $M$, $k$, $D$, $\varphi$, $L$, $S$ and $R$ be as in
\rm
   \bf (\th$^{+}$)\sl.  Assume, moreover, that $T$ is local with 
   nonzero maximal
ideal $M$, $D =L$
   is a proper subfield of $k$, and that $T\subsetneq M^{-1}=(R:M)$.  In
this situation,

  \centerline{$ d_{R} \lneq (d_{D})^{{\varphi}} =
   \left((d_{R})_{{\varphi}}\right)^{\varphi}\,.  $}

  \noindent \vskip 3pt\rm With the notation of Proposition
\ref{prop:3.1}
  (b), take ${\cal{D}} :=\{D\}$, hence ${\cal{R}} =\{R\}$, thus
  $\star_{\cal D}= d_{D}=v_{D}$ and $\star_{\cal R}= d_{R}$.  In this
  situation,
  by Corollary \ref{coro:2.13},
  $(\star_{\cal D})^{\varphi} = v_{R}$.  Therefore, by
  Proposition \ref{prop:3.2} and Example \ref{ex:2.3} (c) and (d),
$(\star_{\cal D})^{\varphi} =
  \left(d_{D}\right)^{\varphi} =
  \left((d_{R})_{{\varphi}}\right)^{\varphi} =v_{R} \gneq
  d_{R}=\star_{\cal R} $.
  \end{exxe}

  Note that it is possible to give an example in which $ \ast \lneq
  \left((\ast)_{\varphi}\right)^{\varphi}$ and $d_R \lneq
  (d_D)^{\varphi}$, even in the case that $D \subsetneq L=k$:

  \begin{exxe} \label{ex:3.4}
      \sl  Let $D$ be a 1-dimensional discrete valuation domain with
quotient
      field $L$.  Set $T:=L[X^2, X^3]$, $M:=X^2L[X]=XL[X] \cap T$
and
      $K:=L(X)$.  Let $\varphi$ and $R$  be as in \rm \bf (\th$^{+}$) \sl
      with $L=k$. Then, $d_R \lneq v_R= ((d_R)_{\varphi})^{\varphi}$.
\rm

      \noindent Since $(d_R)_{\varphi}=d_D=v_D$ and
$(v_D)^{\varphi}=v_R$
      (Corollary \ref{coro:2.13}), we have
$((d_R)_{\varphi})^{\varphi}=v_R$.
      Now consider, for instance, the fractional ideal $T$ of $R$.  We
know,
      from Example \ref{ex:2.17}, that
      $T$ is not a divisorial ideal of $R$, i.e., $T^{d_{R}}= T \lneq
T^{v_{R}}$.
      Thus we have $d_R \lneq v_R= ((d_R)_{\varphi})^{\varphi}$.
  \end{exxe}

  \vskip 4pt

  The next goal is to show that $(t_{R})_{\varphi} = t_{D}$ (but, in
general,
  $t_{R} \lneq (t_{D})^{{\varphi}} =
  \left((t_{R})_{{\varphi}}\right)^{\varphi}$).  We start with a more
general result concerning
  the preservation of the ``finite type'' property.

  \begin{prro}\label{prop:3.4} \sl Let $T$, $K$, $M$,  $k$, $D$,
$\varphi$, $L$, $S$ and
  $R$ be as in \rm \bf (\th$^{+}$)\it.  Let $\ast$ be a given star
  operation on the integral domain $R$.
  \balf \bf \item \it If $\ast$ is
  a star operation of finite type on $R$, then $\ast_{\varphi}$ is a
star
  operation of finite type on $D$.
  \bf \item \it If $\ast$ is any star
  operation on $R$, then $({\ast_f})_{\varphi} =
  (\ast_{\varphi})_{f}\,.$ \ealf \end{prro} \rm

  \noindent {\bf Proof}.
      \bf (a) \rm To prove the statement we will use the following
facts:
      \smallskip

        \hskip 10pt \bf (1) \rm For each integral ideal $I$ of $R$ such that
$M\subset I$,
        $ \left(\frac{I}{M}\right)^{\ast_{\varphi}}=
        \left(\varphi(I)\right)^{\ast_{\varphi}}= \varphi(I^\ast) =
        \frac{I^\ast}{M}\, $ (Proposition \ref{prop:2.7}).  \smallskip

       \hskip 10pt  \bf (2) \rm For each nonzero ideal $I$ of $R$,
$(I+M)^\ast \supseteq I^\ast
       +M$.

      \hskip 10pt  \bf (3) \rm For each nonzero ideal $J$ of $D$ and for
each $y \in L$ with
      $ J \subseteq yD$,   if $F_{y}$ is a finitely generated ideal of
      $R$ such that $F_{y} \subseteq
I_{y}:=\varphi^{-1}\left(y^{-1}J\right)$,
      then $y\varphi(F_{y})$ is a finitely generated ideal of $D$ with
      $y\varphi(F_{y}) \subseteq J$.
      \smallskip

      \smallskip
      For each nonzero ideal $J$ of $D$, we have
      $$
      \begin{array}{rl}
      J^{\ast_{\varphi}} &= \cap

\left\{y\varphi\left(\left(\varphi^{-1}\left(y^{-1}J\right)\right)^{\ast}\right)

  \mid \, y \in L\,,\;  J \subseteq yD \right\} \\ &= \cap
      \left\{y\varphi\left(I_{y}^{\ \ast}\right) \mid \, y \in
      L\,,\; J \subseteq yD \right\} \\
      &= \cap
      \left\{y\varphi\left(\cup \left\{ F_{y}^{\ \ast}\mid  F_{y}
\subseteq
      I_{y},\   F_{y} \in {\boldsymbol{f}}(R) \right\}\right) \mid  y
\in
      L,\  J \subseteq yD \right\} \\
         &= \cap
      \left\{\cup \left\{y\varphi\left( F_{y}^{\ \ast}\right)\mid F_{y}
\subseteq
      I_{y},\   F_{y} \in {\boldsymbol{f}}(R) \right\} \mid y \in L,\ J
      \subseteq yD \right\} \\
       &= \cap
      \left\{\cup \left\{y\frac{ F_{y}^{\ \ast}+M}{M}\mid F_{y}
\subseteq
      I_{y},\   F_{y} \in {\boldsymbol{f}}(R)   \right\} \mid y \in L,\
J \subseteq yD \right\} \\
      &\subseteq \cap \left\{\cup \left\{y\frac{
(F_{y}+M)^{\ast}}{M}\mid F_{y} \subseteq
      I_{y},\   F_{y} \in {\boldsymbol{f}}(R)  \right\} \mid y \in L,\ J
\subseteq yD
      \right\} \\
        &= \cap \left\{\cup
        \left\{y\left(\frac{F_{y}+M}{M}\right)^{\ast_{\varphi}}\mid
F_{y} \subseteq
      I_{y},\   F_{y} \in {\boldsymbol{f}}(R)   \right\} \mid y \in L,\
J \subseteq yD
        \right\} \\
        &= \cap \left\{\cup
        \left\{y\left(\varphi(F_{y})\right)^{\ast_{\varphi}}\mid F_{y}
\subseteq
      I_{y},\   F_{y} \in {\boldsymbol{f}}(R)  \right\} \mid y
        \in L,\ J \subseteq yD \right\}\\
        &= \cap \left\{\cup
        \left\{\left(y\varphi(F_{y})\right)^{\ast_{\varphi}}\mid F_{y}
\subseteq
      I_{y},\   F_{y} \in {\boldsymbol{f}}(R)  \right\} \mid y
        \in L,\ J \subseteq yD \right\}\\
          &\subseteq \cap \left\{\cup \left\{G^{\ast_{\varphi}}\mid G
\subseteq
          J,\ G \in {\boldsymbol{f}}(D) \right\} \mid y \in L,\ J
          \subseteq yD \right\}\\
       &=\cup \left\{ G^{\ast_{\varphi}}\mid G \subseteq J,\ G
       \in {\boldsymbol{f}}(D) \right\} \subseteq J^{\ast_{\varphi}}\,,
      \end{array}
      $$
      where we may assume each $F_y\not\subseteq M$ so that we can use
Fact (1).

      Thus, $J^{\ast_{\varphi}} = \cup \left\{
G^{\ast_{\varphi}}\mid G \subseteq J,\ G
       \in {\boldsymbol{f}}(D) \right\}$.

       \bf (b) \rm Since both $({\ast_f})_{\varphi}$ and
$({\ast_{\varphi}})_f$
       are star operations of finite type on $D$ by (a), it suffices to
show
       that for each nonzero finitely generated ideal $J$ of $D$,
       $J^{({\ast_f})_{\varphi}} = J^{({\ast_{\varphi}})_f}$.  Recall
that if
       $J$ is a nonzero finitely generated ideal of $D$, then
$\varphi^{-1}(J)$
       is a finitely generated ideal of $R$ \cite[Corollary
1.7]{FG:1996}.
       Therefore
         $$
      \begin{array}{rl}
      J^{(\ast_{\varphi})_{f}}= J^{\ast_{\varphi}} &=
      \left\{
                y\varphi\left(
                          \left( \varphi^{-1}
                            \left(y^{-1}J\right) \right)^{\ast} \right)
\mid \, y \in L\,,\; J
                            \subseteq yD \right\} \\
  &=
\left\{y\varphi\left(\left(\varphi^{-1}\left(y^{-1}J\right)\right)^{\ast_{f}}\right)

  \mid \, y \in L\,,\;  J \subseteq yD \right\} \\
  &= J^{(\ast_{f})_{\varphi}}\,.
    \end{array}
      $$
     \vskip -16pt \hfill $\Box$

  \begin{prro}\label{prop:3.5} \it Let $T$, $K$, $M$, $k$, $D$,
$\varphi$,
  $L$, $S$ and $R$ be as in \rm \bf (\th$^{+}$)\it . Then

  \centerline{$ (t_R)_{\varphi}=t_D\,.
  $}
  \end{prro}

  \vskip -6pt \noindent {\bf Proof}.  \rm
  Easy consequence of Corollary \ref{coro:2.13} and Proposition
\ref{prop:3.4} (b).
%
%
%

  \begin{reem} \label{rk:3.6} \sl In the same situation of Example
  \ref{ex:3.4}, choosing $D$ to be a De\-de\-kind domain with infinitely
many
  prime ideals, we have

 \centerline{$ t_{R} \lneq \left(t_{D}\right)^{\varphi}=
  \left((t_{R})_{{\varphi}}\right)^{\varphi}\,.  $} \rm

Using Proposition \ref{prop:3.5} we have $\left(t_{D}\right)^{\varphi}=
  \left((t_{R})_{{\varphi}}\right)^{\varphi}$. We claim that, in the
present
situation, the set of the maximal $t_{R}$--ideals of $R$ coincides with
$\Max(R)$.

Note first that
since dim$(T)=1$, the contraction to $R$ of each nonzero prime ideal of
$T$ has
height 1 \cite[Theorem 1.4]{F},\ so it is a $t_{R}$--prime of $R$
\cite[Corollaire 3, p.  31 ]{J}.  Let $Q\in \Max(R)$.  If
$Q\not\supseteq M$, then $Q$ is the contraction of a prime ideal of $T$,
so $Q$ is a $t_{R}$--prime.  If $Q\supseteq M$, then
$\frac{Q}{M}=(\frac{Q}{M})^{v_D}=\frac{Q^{v_R}}{M}$ by Proposition
\ref{prop:2.7}, and hence we have $Q^{v_R}=Q$.  Therefore, in this case
also, $Q$ is a $t_{R}$--prime.

Note that $M$ is a divisorial prime ideal in $R$, hence in particular
$M$ is a prime
$t_{R}$--ideal and it is contained in infinitely many maximal
($t_{R}$--)ideals, therefore $R$ is not a TV-domain, i.e., $t_R\neq
v_R$\
\cite[Theorem 1.3 and Remark 2.5]{HZ}.  Since
$((d_R)_{\varphi})^{\varphi}=(d_D)^{\varphi}=(v_D)^{\varphi}=v_R$,
automatically we have $((t_R)_{\varphi})^{\varphi}=(t_D)^{\varphi}=v_R$.
Thus, in this example, we have $t_R\lneq (t_D)^{\varphi}$.

  Note also that \sl this example shows that if $\star$ is a star
  operation of finite type on $D$, then $\star^\varphi$ is a star
  operation on $R$, which is not necessarily of finite type \rm  (e.g.
take  $\star := t_{D}= (t_{R})_{\varphi}$).

  \end{reem}

  In the pullback setting that we are considering, it is also natural to
ask
  about the transfer of the property of being a ``stable'' star
operation.

  \begin{prro}\label{prp:3.7} \it Let $T$, $K$, $M$, $k$, $D$,
$\varphi$,
  $L$, $S$ and $R$ be as in \rm \bf (\th$^{+}$) \it and let $\ast$ be a
  star operation on $R$.  Then
%
%

  \centerline{$\tilde{\ast}_{\varphi} =\widetilde{(\ast_{\varphi})}$\,.
}
  \end{prro}

  \noindent {\bf Proof}.  \rm If $D=k$, then since $D=L$ is a
  field, obviously we have $\tilde{\ast}_{\varphi}=
  \widetilde{(\ast_{\varphi})}$. Assume that $D\subsetneq k$.

   Let $J$ be a nonzero integral ideal of $D$ and let
  $I:={\varphi}^{-1}(J)$.  We first show that
  $J^{{\tilde{\ast}}_{\varphi}} \subseteq J^{\widetilde
{(\ast_{\varphi})}}$.\
  By Proposition \ref{prop:2.7},
  $J^{{\tilde{\ast}}_{\varphi}}= \frac{I^{\tilde{\ast}}}{M}$.  Moreover,
  recall that $J^{\widetilde {(\ast_{\varphi})}}=\{y\in D \mid
  yJ_1\subseteq J$ for some finitely generated ideal $J_1$ of $D$ such
  that $J_1^{\ast_{\varphi}}=D\}$ [respectively, $I^{\tilde{\ast}}=
\{x\in
  R \mid xI_1\subseteq I$ for some finitely generated ideal $I_1$ of $R$
  such that $I_1^{\ast}=R\}$].\
  Let $y\in J^{{\tilde{\ast}}_{\varphi}}$.  Then $y=\varphi (x)$ for
  some $x\in I^{\tilde{\ast}}$.  So $xI_1\subseteq I$ for some
  finitely generated ideal $I_1$ of $R$ such that $I_1^{\ast}=R$.
  Set $J_1:=\varphi (I_1)=\frac{I_1+M}{M}$.  Then $J_1$ is nonzero
finitely
  generated, and by Proposition \ref{prop:2.7},
  $J_1^{\ast_{\varphi}}=\frac{(I_1+M)^{\ast}}{M}=\frac{R}{M}=D$.
  Since $xI_1\subseteq I$,  $yJ_1=\varphi (xI_1)\subseteq \varphi
  (I)=J$, and hence $y\in J^{\widetilde {(\ast_{\varphi})}}$.

  Conversely, let $J$ be a nonzero integral ideal of $D$.  If $y\in
J^{\widetilde{(\ast_{\varphi})}}=J^{\widetilde{(\ast_{\varphi})_f}}=J^{\widetilde{(\ast_f)_{\varphi}}}$
  (Proposition \ref{prop:3.4} (b)), then $yJ_1\subseteq J$ for some
finitely generated ideal
  $J_1$ such that $J_1^{(\ast_f)_{\varphi}}=D$.  Set
  $I_1:={\varphi}^{-1}(J_1)$.  Since
  $J_1^{(\ast_f)_{\varphi}}=\frac{I_1^{\ast_{_f}}}{M}=D$ (Proposition
  \ref{prop:2.7}), $I_1^{\ast_f}=R$.  Therefore, there exists a finitely
  generated subideal $I_0$ of $I_1$ such that $I_0^{\ast}=R$.  Write
  $y:=\varphi (x)$ for some $x\in R$.  Since $xI_0\subseteq
xI_1\subseteq
  I:= {\varphi}^{-1}(J)$, $x\in I^{\tilde{\ast}}$, and hence (using
  Proposition \ref{prop:2.7} again) $y\in
  \frac{I^{\tilde{\ast}}}{M}=J^{(\tilde{\ast})_{\varphi}}$.  \hfill
$\Box$

  \begin{coor}\label{cor:3.8} \it Let $T$, $K$, $M$, $k$, $D$,
$\varphi$,
  $L$, $S$ and $R$ be as in \rm \bf (\th$^{+}$)\it.  Then

  \centerline{$ (w_{R})_{\varphi}= w_{D}\,.
  $}
  \end{coor}

  \noindent {\bf Proof}.  \rm Recall that $w_{R} = \widetilde{v_{R}}$
and
  $w_{D} = \widetilde{v_{D}}$.  The conclusion follows from Proposition
  \ref{prp:3.7}.
%
  \hfill
  $\Box$

 \begin{reem} \sl The example considered in Remark \ref{rk:3.6} shows
that we
 can have $w_R\lneq ((w_R)_{\varphi})^{\varphi}=(w_D)^{\varphi}$.  \rm

Since $\Max(R)= \mathcal M(t_{R})$\ (= the set of the maximal
$t_{R}$--ideals, according to the notation in Example \ref{ex:1.3}
(e)),  $w_R={\star}_{\mathcal M(t_{R})}=d_R$.
In particular, $T^{w_R}=T$.  On the other hand, we know that
$((d_R)_{\varphi})^{\varphi}=(d_D)^{\varphi}=(v_D)^{\varphi}=v_R$. Thus
we
have $((w_R)_{\varphi})^{\varphi}=(w_D)^{\varphi}=v_R$.  As we have
already noticed (Example \ref{ex:3.4}), $T$ is not a divisorial ideal of
$R$, i.e., $T^{v_R}\supsetneq T=T^{w_R}$.  Thus, in this case,  we have
$w_R\lneq
(w_D)^{\varphi}$.
\end{reem}

  Since the stable star operation $\tilde{\ast}$ is a particular type of
  spectral star operation, the next goal is a possible extension of
  Proposition \ref{prp:3.7} to the case of spectral star operations.  We
  start with the following:

  \begin{leem} \label{lm:3.9} Let $T$, $K$, $M$, $k$, $D$, $\varphi$,
  $L$, $S$ and $R$ be as in \rm \bf (\th$^{+}$)\it.  Assume that
  $D\subsetneq k$.  \balf\bf \item \it Let $P$ be a prime ideal of $R$
  containing $M$.  Set $Q: ={\varphi}(P)$ and $R_{(P, \varphi)}:
  ={\varphi}^{-1}(D_Q)$.  Then $R_{(P, \varphi)}=R_P\cap T$.  \bf \item
  \it Let $\Delta(\neq\emptyset) \subseteq \mbox{\rm Spec}(R)$ and
assume
  that $\ast :={\star}_{\Delta}\in \mbox{\bf Star}(R)$.  Set
  ${\Delta}_1:=\{P\in\Delta\mid P\supseteq M\}$.  For each nonzero
  integral ideal $I$ of $R$ containing $M$, we have

      \vskip -5pt \centerline{$I^\ast=\cap \{IR_{(P, \varphi)} \mid P\in
{\Delta}_1\}\,.  $}

      (Note that ${\Delta}_1\neq\emptyset$.)  \ealf
  \end{leem}
  {\bf Proof.} \bf (a) \rm is straightforward.

  \bf (b) \rm If $M=(0)$, then $\Delta={\Delta}_1$ and
  $R_{(P,\varphi)}=R_P$, so it trivially holds. Assume that $M\neq
  (0)$.
  Let $I$ be an integral ideal of $R$ containing
  $M$. Recall that for each $P\in \Delta \setminus {\Delta}_1$,
  there exists a unique $P' \in \Spec(T)$ such that $P' \cap R=P$
  and $R_P=T_{P'}$ \cite[Theorem 1.4]{F}, hence in particular
  ${\Delta}_1\neq\emptyset$ (otherwise ${\star}_{\Delta}$ would not
  be a star operation on $R$).  We have $$
  \begin{array}{rl}
      I^\ast &=\cap \{IR_P \mid P\in\Delta \}=
      (\cap \{IR_P \mid P\in {\Delta}_1 \}) \cap (\cap \{IR_P \mid
P\in\Delta\setminus
      {\Delta}_1\}) \\
      &=(\cap \{IR_P \mid P\in{\Delta}_1 \}) \cap (\cap \{R_P \mid
P\in\Delta\setminus
      {\Delta}_1\})
     \\
    &\supseteq (\cap \{IR_P \mid P\in{\Delta}_1 \}) \cap T \supseteq
    \cap \{IR_{(P, \varphi)} \mid  P\in{\Delta}_1 \}\,.
    \end{array}
    $$

  Conversely, let $x\in I^\ast $ and let $P\in {\Delta}_1$ (which is
  nonempty). Then there exists $s\in R\setminus P$ such that $sx\in
  I$.  Since ${\varphi}(s)\in D\setminus {\varphi}(P)$,
  ${\varphi}(s)$ is a unit element of $D_{\varphi(P)}$, and hence there
exists
  $t\in R_{(P, \varphi)} $ such that $\varphi(t)\varphi(s)=1$, or
  equivalently, $ts-1\in M$.  Put $ts-1=:m\in M$, then
$tsx=(1+m)x=x+mx$.
  Since $tsx\in IR_{(P, \varphi)} $ and $mx\in MI^\ast \subseteq MR=
  M\subseteq I \subseteq IR_{(P, \varphi)} $,  we have $x= tsx -mx\in
IR_{(P,
  \varphi)}$.  \hfill $\Box$

  \begin{prro}\label{prop:3.10} \it Let $T$, $K$, $M$, $k$, $D$,
$\varphi$,
  $L$, $S$ and $R$ be as in \rm \bf (\th$^{+}$)\it.  Let $\Delta$ be a
  nonempty set of prime ideals of $R$ and assume that $\ast:
  ={\star}_{\Delta}\in \mbox{\bf Star}(R)$.  Set
  ${\Delta}_{\varphi}:=\{{\varphi}(P)\mid P\in\Delta,\, P\supseteq M\}\
  (\subseteq \mbox{\rm Spec}(D))$.  Then

\centerline{$({{\star}_{\Delta}})_{\varphi}={\star}_{{\Delta}_{\varphi}}\,.
  $}
  \end{prro}

  {\bf Proof.} \rm If $D=k$, then since $D=L$ is a field, we obviously
have
  $({\star}_{\Delta})_{\varphi}={\star}_{{\Delta}_{\varphi}}$.  Assume
  that $D\subsetneq k$, then ${\Delta}_{\varphi}\neq\emptyset$.  Let $J$
  be a nonzero integral ideal of $D$ and let $I:={\varphi}^{-1}(J)$.
Set
  ${\Delta}_1=\{P\in\Delta \mid P\supseteq M\}$, hence
  ${\Delta}_{\varphi}=\{\varphi(P) \mid P\in {\Delta}_1\}$.  Since $I$
is
  an integral ideal of $R$ containing $M$, $I^\ast =\cap \{IR_{(P,
  \varphi)} \mid P\in{\Delta}_1 \}$ by Lemma \ref{lm:3.9} (b), and so,
  using Proposition \ref{prop:2.7}, we have
  $J^{\ast_{\varphi}}=\varphi(I^\ast)= \cap \{ \varphi(I)D_{\varphi(P)}
  \mid P \in {\Delta}_1 \} = \cap \{JD_{\varphi(P)} \mid P\in{\Delta}_1
  \}= \cap \{JD_{Q} \mid Q\in \Delta_{\varphi} \}=J^{
  \star_{\Delta_{\varphi}} }$.  \hfill $\Box$

  \vskip 10pt

  \begin{reem}\label{rk:3.11} \bf (1) \rm Note that from Proposition
  \ref{prop:3.10} we can deduce another proof of Proposition
  \ref{prp:3.7}.  As a matter of fact, for each star operation
  $\ast$ on $R$, $\tilde{\ast} = {\star}_{\Delta}$, where
  $\Delta:={\cal M}(\ast_{f})$ (Example \ref{ex:1.3} (e)).
In the present situation, ${\Delta}_1:=\{ P \in {\cal M}(\ast_{f})
  \mid P \supseteq M \}$.  By using Proposition \ref{prop:2.7} and
  Proposition \ref{prop:3.4}  (b), it is easy to see that

  \centerline{$ P \in \Delta_{1}\;\;\; \Leftrightarrow \;\;\;
Q:=\varphi(P) \in
  {\cal
  M}((\ast_{\varphi})_{f})\,.$}

  \bf (2) \sl Note that if $\star:= {\star}_{\Delta}$ is a spectral star
  operation on $D$, then $\star^{\varphi}$ is not necessarily a spectral
star
  operation on $R$ \rm (in particular, $({\star}_{\Delta})^{\varphi}
\neq
  \star_{\Delta^{\varphi}}$, where $\Delta^{\varphi} := \{P \in
  \mbox{Spec}(R) \mid \varphi(P) \in \Delta \}$).

  To show this fact, let $D$ be a 1-dimensional discrete valuation
  domain with quotient field $L$ and maximal ideal $N$. Let
  $T:=L[\![X^2, X^3]\!]$ and let $M:=X^2L[\![X]\!]=XL[\![X]\!]\cap T$.  Under these hypotheses, let
  $R$ be the integral domain  defined (as a pullback of type \bf
  (\th$^{+}$)\rm ) from $D$, $T$ and the canonical projection
  $\varphi: T \rightarrow L$.  Then, $R$ is a 2-dimensional
  non-Noetherian local domain.  Let $\Delta: =\mbox{Max}(D)= \{N\}$.
  Then $\star:={\star}_{\Delta}=d_D=v_D$ and
  ${\star}^{\varphi}=(v_D)^{\varphi}=v_R$ (Corollary
  \ref{coro:2.13}).  Since ${\Delta}^{\varphi}=\mbox{Max}(R)$,
  ${\star}_{{\Delta}^{\varphi}}=d_R$.  Suppose that
  ${\star}^{\varphi}$ is spectral, then by Proposition
  \ref{prop:3.10} and Proposition \ref{prop:2.9}, we have
  necessarily that ${\star}^{\varphi}$ coincides with
  ${\star}_{{\Delta}^{\varphi}}$, i.e., $v_R= {\star}^{\varphi}=
  {\star}_{{\Delta}^{\varphi}}= d_R$.  This is a contradiction,
  since $T^{v_R}=(R:_{K}(R:_{K}T))=(R:_{K}M)\supseteq
  L[\![X]\!]\supsetneq T=T^{d_R}$.

  \end{reem}

  \begin{prro}\label{prop:3.12} \it Let $T$, $K$, $M$, $k$, $D$,
$\varphi$,
  $L$, $S$ and $R$ be as in \rm \bf (\th$^{+}$)\it .  If $\ast$ is an
a.b.
  [respectively, e.a.b.] star operation on $R$, then $\ast_{\varphi}$ is
an a.b.
  [respectively, e.a.b.]
  star operation on $D$.
  \end{prro}

  \vskip -4pt {\bf Proof.} Let $J$ be a nonzero finitely generated ideal
of $D$
  and let $J_1$, $J_2$ be two arbitrary nonzero ideals of $D$ such that
  $(JJ_1)^{\ast_{\varphi}}\subseteq (JJ_2)^{\ast_{\varphi}}$.  Set
  $I:={\varphi}^{-1}(J)$, $I_i:={\varphi}^{-1}(J_i)$ for $i=1, 2$.
Since
  $J$ is finitely generated and $IS=S$ (because $I \supset M$ and $M$ is
a maximal ideal of $S$), there exists a finitely generated subideal
$I_0$ of $I$ such that ${\varphi}(I_0)=J$ and $I_0S=S$.  Then, by
Proposition \ref{prop:2.7},
we have $(I_0I_1+M)^\ast\subseteq (I_0I_2+M)^\ast$.
  Note that $I_0I_i\supseteq I_0M=I_0MS=I_0SM=SM=M$ for $i=1,2$, thus
we have $(I_0I_1)^\ast\subseteq (I_0I_2)^\ast$.  Since $I_0$ is finitely
  generated and $\ast$ is an a.b. star operation,
  ${I_1}^\ast\subseteq {I_2}^\ast$ and so
${J_1}^{\ast_{\varphi}}\subseteq
  {J_2}^{\ast_{\varphi}}$.
  The statement for the e.a.b. case follows from Proposition
  \ref{prop:3.4} (b) and from the fact that $\ast$ is e.a.b. if and only
if
  $\ast_f$ is a.b..\hfill $\Box$

  \begin{reem} \label{rk:3.13} \bf (1) \rm
%
%
%
  \sl Under the assumption of Proposition \ref{prop:3.12}, if $v_R$ is
  e.a.b., then $(v_R)_{\varphi}=v_D$ is e.a.b..  In other words, if $R$
is
  a $v_{R}$--domain, then $D$ is a $v_{D}$--domain \rm \cite[page
418]{G}.

  \bf (2) \sl Let $\star$ be an a.b. [respectively, e.a.b.] star
operation on
  $D$.  Then, in general, ${\star}^{\varphi}$ is not an a.b.
  [respectively, e.a.b.] star operation on $R$.  \rm

  To show this fact, take $D$, $T$ and $R$ as in Remark
  \ref{rk:3.11} (2).  Since $D$ is a 1-dimensional discrete
  valuation domain, its unique star operation $d_D(= b_{D}=v_D)$ is
  an a.b. star operation (and hence an e.a.b. star operation). Since
  $R$ is not integrally closed (because $X\in K\setminus R$ is
  integral over $R$), $R$ has no e.a.b. star operations (and hence
  no a.b. star operations).

  Note that it is possible to give an example of this phenomenon also
with
  $R$ integrally closed.

  \end{reem}

  \begin{exxe} \label{ex:3.15} \sl Let $D$ be a 1-dimensional discrete
  valuation domain with quotient field $L$, let $T:=L[X, Y]$ and
  $M:=(X, Y)L[X, Y]$.  Under these hypotheses, let $R:= D+(X, Y)L[X,
  Y]$ be the integral domain defined (as a pullback of type \bf
  (\th$^{+}$)\sl ) from $D$, $T$ and the canonical projection $\varphi:
T
  \rightarrow L$.  Then $(b_D)^{\varphi}$ is not e.a.b. (and hence not
  a.b.) on $R$.  \rm

  Note that $M$ is a
  divisorial ideal of $R$ of finite type, in fact, $M=I^{v_R}$, where
  $I:=(X, Y)R$.  Now, choose $a_1, a_2\in D\setminus (0)$ such that
  $a_1D\not\subseteq a_2D$ (e.g. put $a_1:=a$, $a_2:=a^2$, where $a$ is
a
  nonzero nonunit element in $D$).  Set $I_1:=a_1R$ and $I_2:=a_2R$.
Then
  $(II_i)^{v_R}=(a_iI)^{v_R}=a_iI^{v_R}=a_iM=M$ (where the last equality
  holds because $a_i$ is a unit in $T$) for each $i=1,2$.  Thus we have
  $(II_1)^{v_R}= (II_2)^{v_R}$.  On the other hand, since
  $(I_i)^{v_R}=I_i= a_iR=a_i(D+M)=a_iD+M$ for each $i=1,2$, and
  $a_1D\not\subseteq a_2D$, we have that $(I_1)^{v_R}\not\subseteq
  (I_2)^{v_R}$.  Therefore, $v_R$ is not an e.a.b. operation.  Since $D$
  is a 1-dimensional discrete valuation domain, the unique star
operation
  $d_D=b_{D}= v_D$ on $D$ is an a.b. star operation (and hence an e.a.b.
  star operation), but $v_R=(v_D)^{\varphi}$ (Corollary \ref{coro:2.13})
  is not e.a.b. (and hence not a.b.).
  \end{exxe}

  Recall that given an integral domain $T$, \it the paravaluation
  subrings of $T$\rm , in Bourbaki's sense \cite[Chap.  6, \S 1,
  Exercise 8]{B}, are the subrings of $T$ obtained as an
  intersection of $T$ with a valuation domain having the same quotient
field
  as $T$.  It is easy to see that if $R$ is a subring of $T$ then the
  integral closure of $R$ in $T$ coincides with the intersection of all
  the paravaluation subrings of $T$ containing $R$ \cite[Chap.  6, \S 1,
  Exercise 9]{B}.

  \begin{leem} \label{le:3.15} \it Let $T$, $K$, $M$, $k$, $D$,
$\varphi$,
  $L$, $S$ and $R$ be as in \rm \bf (\th$^{+}$)\it .  Assume that $D
  \subsetneq L=k$.  Assume, moreover, that $D$ is integrally closed (or
  equivalently, that $R$ is integrally closed in $T$).  Let
  ${\cal{P}}:={\cal{P}}(R, T)$\ [respectively, ${\cal{V}}$;\
  ${{\cal{V}}_{1}}$;\ ${\cal{W}}$]\ be the set of all the paravaluation
  subrings of $T$ containing $R$\ [respectively, the set of all
valuation
  overrings of $R$;\ the set of all valuation overrings $(V_{1}, N_{1})$
  of $R$ such that $N_{1} \cap R \supseteq M $;\ the set of all the
  valuation overrings of $D$].\ Set $b_{R, T} := \star_{{{\cal{P}}}}$\
  [respectively, $b_{R} := \star_{{\cal{V}}}$;\ $\star_{1} := \star_{
  {{\cal{V}}}_{1} }$;\ $b_{D} := \star_{{\cal{W}}} $ ].\ Then

  \balf \bf \item \it $b_{R, T}$ [respectively, $b_{D}$] is a star
operation
  on $R$ [respectively, on $D$];\ $b_{R} $ and $\star_1$ are semistar
  operations on $R$.  Moreover,

  \centerline{$ b_{R, T} \leq \star_1 \wedge \star_{\{T\}} \leq b_{R}\,.
  $}

   \bf \item \it $(b_{R,T})_{\varphi}=b_D$.

  \bf \item \it If $R$ is integrally closed (which happens if $T$ is
  integrally closed), then $\star_1 \wedge \star_{\{T\}}$ and $b_{R}$
are
  star operations on $R$.\ Moreover, $(b_R)_{\varphi}=b_D$ and $b_R \leq
  (b_D)^{\varphi}$.

  \bf \item \it If $T:=V$ is a valuation domain, then $b_{R, T}
=\star_1=
  \star_1 \wedge \star_{\{T\}}=b_{R}$ .  \ealf

  \end{leem}

  {\bf Proof.} Note that if $(V_{2}, N_{2}) \in {{\cal{V}}}\setminus
  {{\cal{V}}}_{1}$, then $N_{2} \cap R \not\supseteq M $, and so there
  exists a unique prime ideal $Q_{2}$ in $T$ such that $R _{N_{2} \cap
R}
  = T_{Q_{2}}$ \cite[Theorem 1.4]{F}.  Therefore $V_{2}\supseteq R
_{N_{2}
  \cap R} = T_{Q_{2}} \supseteq T$.

  \bf (a) \rm The first part of this statement is an obvious
  consequence of the definitions and the assumption that
  $R$ is integrally closed in $T$ (and equivalently, $D$ is
  integrally closed \cite[Corollary 1.5]{F}).  For each $I \in
  \overline{\boldsymbol{F}}(R)$, we have $I^{b_{R}} = \cap \{ IV \mid V \in
  {\cal{V}}\} = \left(\cap \{ IV_{1} \mid V_{1} \in
  {\cal{V}}_{1}\}\right) \cap \left(\cap \{ IV_{2} \mid V_{2} \in
  {\cal{V}} \setminus {\cal{V}}_{1}\}\right) \supseteq \left(\cap \{
  IV_{1} \mid V_{1} \in {\cal{V}}_{1}\}\right) \cap IT =
  I^{\star_{1}} \cap I^{\star_{\{T\}}} \supseteq \left(\cap \{
  I(V_{1}\cap T) \mid V_{1} \in {\cal{V}}_{1}\}\right) \supseteq
  (\cap \{ I(V \cap T)$
   $\mid V \in {\cal{V}}\}) = I^{b_{R, T}}$.

  \bf (b) \rm Note that since $L$ is a field, the paravaluation
subrings
  of $L$ containing $D$ coincide with the valuation rings in $L$
  containing $D$ \cite[Chap.  6, \S 1, Exercise 8 d)]{B}.  Moreover, if
  $W$ is a valuation overring of $D$, then $\varphi^{-1}(W)$ is a
  paravaluation subring of $T$ containing $R$ \cite[Chap.  6, \S 1,
  Exercise 8 c)]{B}.  On the other hand, if $V'\cap T$ is a
paravaluation
  subring of $T$ (where $V'$ is a valuation domain in the field $K$,
  quotient field of $R$), then necessarily $\varphi(V'\cap T)$ is a
  paravaluation subring of $\varphi(T) =L$, i.e., it is a valuation
domain
  in $L$ containing $D$ \cite[Chap.  6, \S 1, Exercise 8 d)]{B}.
  Therefore, for each $J \in \boldsymbol{F}(D)$,
$\varphi^{-1}(J^{b_{D}})
  = (\varphi^{-1}(J))^{b_{{R,T}}}$.  Now, we can conclude, since
  we know that for each $J \in \boldsymbol{F}(D)$,
  $J^{(b_{{R,T}})_{\varphi}} =
  ({(\varphi^{-1}(J))^{b_{{R,T}}}})/{M}\,$ (Proposition \ref{prop:2.7}).

    \bf (c) \rm   If $R$ is
  integrally closed, then $b_{R}$ is a star operation on $R$
  \cite[Corollary 32.8]{G}, and so by (a) it follows that  $\star_1
\wedge
  \star_{\{T\}}$ is also a star operation on $R$.

  Let $\mathcal W=\{W_{\lambda} \,|\,
\lambda\in\Lambda\}$. For each $\lambda\in\Lambda$, let
$R_{\lambda}:={\varphi}^{-1}(W_{\lambda})$. Then, by the argument
used in the proof of (b), we have $\mathcal
P=\{R_{\lambda}\,|\,\lambda\in\Lambda\}$. Denote by $A'$ the
integral closure of an integral domain $A$. Since $R_{\lambda}$ is
integrally closed in $T$, $R_{\lambda}={R_{\lambda}}' \cap T$. Let
$\iota'_{\lambda}: R_{\lambda} \hookrightarrow {R_{\lambda}}'$ and
$\iota_{\lambda}: R_{\lambda} \hookrightarrow T$ be the canonical
embeddings, and set
$\ast_{\lambda}:=(b_{{R_{\lambda}}'})^{\iota'_{\lambda}}\wedge
(d_T)^{\iota_{\lambda}}$ for each $\lambda \in \Lambda$\ (note
that $({d_T})^{\iota_{\lambda}}$ coincides with the semistar
operation ${\star}_{\{T\}}$ on ${R_{\lambda}}$ ).\ Then
$\ast_{\lambda}$ is a star operation on $R_{\lambda}$ (see also
\cite[Theorem 2]{Anderson:1988}).

{\bf Claim 1.} \sl  Let $I$ be an integral ideal of $R$ properly
containing $M$.  Then $(IR_{\lambda})^{\ast_{\lambda}}=IR_{\lambda}$.\rm

Let $Q_{\lambda}$ be the maximal ideal of the valuation domain
$W_{\lambda}$. If
${\varphi}(IR_{\lambda})=\frac{IR_{\lambda}}{M}\neq yQ_{\lambda}$
for all $y\in L\setminus (0)$,  then since
${\varphi}(IR_{\lambda})$ is a divisorial ideal of the valuation
domain $W_{\lambda}$, ${\varphi}(IR_{\lambda}) =
{\varphi}\left((IR_{\lambda})^{\ast_{\lambda}}\right)$ (and
hence, $(IR_{\lambda})^{\ast_{\lambda}}=IR_{\lambda}$) by
Proposition \ref{prop:2.7}.  \ Assume that
${\varphi}(IR_{\lambda})=\frac{IR_{\lambda}}{M}=yQ_{\lambda}$ for
some $y\in L\setminus (0)$.  Choose $s\in S\setminus M$ such that
${\varphi}(s) =y$ and let
$P_{\lambda}:={\varphi}^{-1}(Q_{\lambda}) \subsetneq R_{\lambda}$.
Then $IR_{\lambda}=sP_{\lambda}+M$, and by the Claim in the proof
of Proposition \ref{prop:2.7}, we have
$(IR_{\lambda})^{\ast_{\lambda}}=s{P_{\lambda}}^{\ast_{\lambda}}+M$.
By (b),  $R_{\lambda}= V_{\lambda} \cap T$ for some valuation
overring $V_{\lambda}$ of $R$, which has center $P_{\lambda}$ on
$R_{\lambda}$, thus
${P_{\lambda}}^{\ast_{\lambda}}=(P_{\lambda}{R_{\lambda}}')^{b_{{R_{\lambda}}'}}
\cap P_{\lambda}T\subseteq P_{\lambda}V_{\lambda}\cap
T=P_{\lambda}$. Therefore, in either case, we have
$(IR_{\lambda})^{\ast_{\lambda}}=IR_{\lambda}$.

{\bf Claim 2.} $(b_R)_{\varphi}\leq (b_{R, T})_{\varphi}$\ ($=b_D$
by (b)).

It suffices to show that for each nonzero integral ideal $J$ of
$D$, $J^{(b_R)_{\varphi}}\subseteq J^{(b_{R, T})_{\varphi}}$,
i.e., for each integral ideal $I$ of $R$ properly containing $M$,
$I^{b_R}\subseteq I^{b_{R, T}}$. Let $I$ be such an ideal. Then
$I^{b_{R,T}}=\cap \{IR_{\lambda}\,|\, \lambda\in\Lambda\}=\cap
\{(IR_{\lambda})^{*_{\lambda}}\,|\, \lambda\in\Lambda\}=\cap
\{(I{R_{\lambda}}')^{b_{{R_{\lambda}}'}}\cap IT \,|\,
\lambda\in\Lambda\}=\cap
\{((I{R_{\lambda}}')^{b_{{R_{\lambda}}'}}\cap T \,|\,
\lambda\in\Lambda\}=\cap \{((I{R_{\lambda}}')^{b_{{R_{\lambda}}'}}
\,|\, \lambda\in\Lambda\}\cap T= \cap \ \{\cap
 \{IV\,|\, V \in {\mathcal V}_{\lambda}:= \{\mbox{valuation
overrings of } {R_{\lambda}}'\} \ \} \, \mid $ $
\;\lambda\in\Lambda\}\cap T
\supseteq \cap \{IV\,|\, V\in \mathcal V\}=I^{b_R}$.

 Therefore, by Claim 2, (a), and the first part of (c), we conclude that
 $(b_R)_{\varphi}=b_D$.  Finally, by  Theorem
  \ref{thr:2.11},  we have
  $b_{R} \leq ((b_{R})_{\varphi})^\varphi =(b_{D})^\varphi $.

  \bf (d) \rm If $T:=V$ is a valuation domain, then each valuation
  overring of $R$ is comparable with $V$.  As a matter of fact, if $V'$
is
  a valuation overring of $R$ and $V'\not\subseteq V$, then there exists
  $y \in V'\setminus V$, hence $y^{-1} \in M$, thus for each $v \in V$,
we
  have $v = v(y^{-1}y) = (vy^{-1})y\in MV'\subseteq V'$.  Therefore, $V
  \subseteq V'$.  From this observation, we immediately deduce that
when
  $T$ is a valuation domain, $b_{R, T} = \star_1 = \star_1 \wedge
  \star_{\{T\}} =b_{R}$.  \hfill $\Box$

  \begin{reem} \rm
  In a pullback situation of type \bf
(\th$^{+}$)\rm ,
when $D$ is integrally closed, we have already
noticed that if $R$ is not integrally closed, then there is no
hope that $(b_{D})^\varphi = b_{R}$ (Remark \ref{rk:3.13} (2)).
More explicitly, Example \ref{ex:3.15}
shows that we can have $b_R\lneq (b_D)^{\varphi}$, even when $R$ is
integrally closed.   The next example shows that $b_R\lneq
(b_D)^{\varphi}$
is possible even under the hypotheses of Lemma \ref{le:3.15} (d).
  \end{reem}

\begin{exxe} \sl Let $T:=V$ be a valuation domain with
maximal ideal $M$ and let $\varphi: V \rightarrow V/M =: k$ be the
canonical projection. Let $D$ be a Dedekind domain with infinitely
many prime ideals and with quotient field $L=k$.  Set $R : =
\varphi^{-1}(D)$. Then $b_R \lneq
(b_D)^{\varphi}$. \rm

By the same
argument as in Remark \ref{rk:3.6}, we can see that $R$ is not a
TV-domain,
i.e., $t_R\neq v_R$. Meanwhile, since $R$ is a Pr\"{u}fer domain,
$b_R=d_R=t_R$, and since $D$ is a Dedekind domain, $b_D=v_D$ and so
$(b_D)^{\varphi}=(v_D)^{\varphi}=v_R$ (Corollary \ref
{coro:2.13}). Therefore, we have $b_R=t_{R}\lneq
v_{R}=(b_D)^{\varphi}$.
\end{exxe}

  \vspace{.1in}

  \vspace{.1in}
  \vspace{.1in}





   \end{document}